%% file: main.tex
\documentclass[a4paper,11pt]{article}

\input{packages}%
\input{commands}%

\graphicspath{{figs}}

\title{Multiscale Segmentation via \\Bregman Distances and Nonlinear Spectral Analysis}
\author[1,2]{Leonie Zeune\footnote{\href{mailto:l.l.zeune@utwente.nl}{l.l.zeune@utwente.nl}}}
\author[2]{Guus van Dalum}
\author[2]{Leon W.M.M. Terstappen}
\author[1]{S.A. van Gils}
\author[1]{Christoph Brune}
\affil[1]{Department of Applied Mathematics}
\affil[2]{Department of Medical Cell BioPhysics}
\affil[ ]{MIRA Institute for Biomedical Technology and Technical Medicine}
\affil[ ]{University of Twente, The Netherlands}
\date{\today}

\begin{document}

\maketitle

\begin{abstract} In biomedical imaging reliable segmentation of objects (e.g. from small cells up to large organs) is of fundamental importance for automated medical diagnosis. New approaches for multi-scale segmentation can considerably improve performance in case of natural variations in intensity, size and shape. This paper aims at segmenting objects of interest based on shape contours and automatically finding multiple objects with different scales. The overall strategy of this work is to combine nonlinear segmentation with scales spaces and spectral decompositions. We generalize a variational segmentation model based on total variation using Bregman distances to construct an inverse scale space. This offers the new model to be accomplished by a scale analysis approach based on a spectral decomposition of the total variation. As a result we obtain a very efficient, (nearly) parameter-free multiscale segmentation method that comes with an adaptive regularization parameter choice. To address the variety of shapes and scales present in biomedical imaging we analyze synthetic cases clarifying the role of scale and the relationship of Wulff shapes and eigenfunctions. To underline the potential of our approach and to show its wide applicability we address three different experimental biomedical applications. In particular, we demonstrate the added benefit for identifying and classifying circulating tumor cells and present interesting results for network analysis in retina imaging. Due to the nature of nonlinear diffusion underlying, the mathematical concepts in this work offer promising extensions to nonlocal classification problems.\\

\textbf{Keywords}: 
 Multiscale Segmentation, Chan-Vese Method, Bregman Iteration, Total Variation, Inverse Scale Space, Nonlinear Spectral Methods, Eigenfunctions, Wulff Shapes, Circulating Tumor Cells
\end{abstract}

%
%
%
%
\section{Introduction}

In mathematical imaging the problem of segmentation refers to the process of automatically detecting regions, objects or patterns of interest. This is of particular importance in biomedical imaging for cell or organ quantification as well as in materials science and engineering. The main goal of this work is to present a new multiscale segmentation method combining variational inverse scale-space methods with nonlinear spectral analysis.

\textit{Image segmentation.} The task of image segmentation can sometimes be addressed by simple methods which are solely based on intensity or histogram thresholds. However those methods quickly fail in more complex experimental scenarios, where challenges in terms of different sizes, intensities, contrasts or uncertainties like noise occur. Thus, a class of commonly used mathematical techniques to overcome some of the difficulties are nonlinear variational methods. In general there are two different ways of describing a region, either by its edge or by the interior of the region. Therefore two main concepts of addressing segmentation by an energy minimization problem evolved in the previous decades: edge-based and region-based segmentation. While edge-based segmentation separates regions based on discontinuity information \cite{Kass1988,Kichenassamy1995,Caselles1997}, region-based segmentation separates by similarity measures within the regions. In this paper we focus on a region-based approach.

\begin{figure}[t]
  \centering 
  \subfigure[Sizes]{\includegraphics[height=0.15\textwidth]{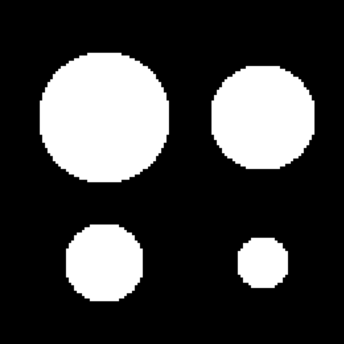}}\qquad
  \subfigure[Intensities]{\includegraphics[height=0.15\textwidth]{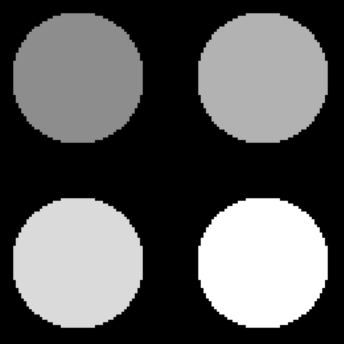}}\qquad
  \subfigure[Shapes]{\includegraphics[height=0.15\textwidth]{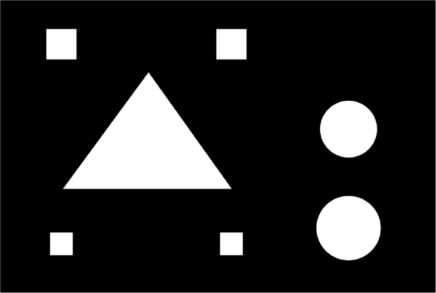}}\\
  \caption{Minimal examples where multiscale image segmentation could offer added value}
\label{fig:overview_goal_segmentation}
\end{figure}

The proposed method in this work is mainly built upon the Chan-Vese (CV) model \cite{Chan2001}. For a given image function $f:\Omega \rightarrow \R$ the domain $\Omega \subset \R^d$ should be separated into two regions $\Omega_1$ and $\Omega_2 := \Omega \setminus \Omega_1$ via the following variational energy minimization method
\vspace{-2pt}
\begin{equation}\label{eq:CVorig}
\int_{\Omega_1} (f(x) - c_1)^2 \mathrm{d}x + \int_{\Omega_2} (f(x) - c_2)^2 \mathrm{d}x + \alpha \cdot Per(C)\longrightarrow \min_{C,c_1,c_2}
\end{equation}
\vspace{-2pt}
where $C$ is the desired contour separating $\Omega_1$ and $\Omega_2$ and $c_i$ stands for a desired average intensity value within region $\Omega_i$. $Per(C)$ denotes the length of the contour $C$ separating $\Omega_1$ and $\Omega_2$ and can be controlled via a regularization parameter $\alpha > 0$. An easy example is the image $f$ given in Figure \ref{fig:overview_goal_segmentation}(a) where the domain $\Omega$ can be segmented with respect to foreground and background. In contrast to simple histogram thresholding methods this holds true even in the case of strong noise if the parameter $\alpha$ is chosen adequately. A generalization of this model was introduced in \cite{Vese2002} and in \cite{Chan2000} the model has also been extended to vector-valued images like color images. Recently in \cite{zosso2015} an extension of the CV model that can handle artifacts and illumination bias in images has been proposed.

\textit{Challenges/Questions.} The CV model is very useful in segmenting regions of interest which have very similar intensity values, e.g. Figure \ref{fig:overview_goal_segmentation}(a). However, automatically detecting single objects based on their size is more challenging. Even with a varying parameter $\alpha$ controlling the contour length (forward scale-space), it is for example not possible to detect the smallest object as a singleton. A similar challenge occurs when segmenting separate objects due to their intensity values, e.g. Figure \ref{fig:overview_goal_segmentation}(b). Increasing the number of constants $c_i$ to four is suboptimal because we usually a-priori don't know the number of objects. Varying a threshold or parameter $\alpha$ could lead to a correct segmentation but the estimated intensity constants $c_1$, $c_2$ will likely be incorrect. \textit{Hence, is it possible to automatically detect multiple scales in a nonlinear variational image segmentation model, for instance with respect to different object sizes or object intensities? Can the segmentation of an image automatically be decomposed with respect to those scales?}

Many region-based segmentation methods only use constraints on the contour length or curvature as regularization. However, in view of shape optimization and dictionary learning an approach that could also automatically separate objects with respect to their shape (cf. Figure \ref{fig:overview_goal_segmentation}(c)) would be very interesting. \textit{Hence, what is the role of geometric shapes in a multiscale segmentation approach?}

\textit{Scale-space methods.}
In the previous decades there has been a continuous interest in the analysis of different scales and the construction of scale spaces in imaging. In general it is desired to automatically detect all scales present in an image and simultaneously determine which scales are informative and contribute most to the image. For segmentation this problem is addressed in the fundamental works by Witkin and Koenderink \cite{Witkin1983,Koenderink1984}. In relation to those, several methods to detect and analyze interesting scales have been proposed, see for example \cite{Lifshitz1990,Vincken1997,Tabb1997,Lindeberg1998,Florack2000,Letteboer2004}. The underlying scale-space that is examined is defined by a linear diffusion process. A drawback of those approaches is that linear diffusion smoothes edge information and is therefore in general not suitable for applications where one is interested in retaining sharp edge information. Especially in biomedical image applications this is often the case. Therefore those theories were extended to non-linear diffusion processes, see \cite{Niessen1997,Niessen1999,Dam2000}. A drawback of these approaches is that their analysis of scales is not fully automatic and can only be used in a forward approach, thus going from fine to coarse scales and then trying to find a backward relation. In this work we concentrate on nonlinear diffusion processes for segmentation where scale automatically relates to intensity and size.

A prominent example of a variational method for nonlinear diffusion is the ROF model \cite{Osher2005}. With increasing regularization parameter $\alpha$ a sequence of functionals generates a nonlinear forward scale space flow that filters signals from fine to coarse. However in this process the total variation regularization functional is known to lead to a systematic contrast loss in the filtered image $u$ \cite{Meyer2001}, whereas the main discontinuities in the signal remain at their position in the domain. To tackle the problem of intensity loss Osher et al. proposed in \cite{Osher2005} an iterative contrast enhancement procedure based on Bregman distances. This approach is known to generate a nonlinear inverse scale space flow generating filtered signals from coarse to fine and with improved quality. This idea was successfully applied to more general inverse problems. \cite{Burger2007,Brune2011,Benning2013}

\textit{Spectral methods.} Recently Gilboa \cite{Gilboa2013,Gilboa2014} developed a framework to detect scales based on the nonlinear total variation diffusion process. The total variation is known to retain edge information while smoothing the signal apart from the edges. In this framework scales are detected based on a spectral decomposition of the given image into TV eigenfunctions \cite{Meyer2001}. This concept does not only hold true for higher-order regularization functionals \cite{Papafitsoros2015511,Poeschl2015,Benning2013} but more generally for convex, one-homogenous functionals $J$ with corresponding nonlinear eigenvalue equations \cite{Burger2015,Burger2016} of the following form
\vspace{-5pt}
\begin{equation*}
	\lambda u \in \partial J(u) \ 
\end{equation*}
where $\partial J(\cdot)$ denotes the subdifferential of $J$. When minimizing for instance the total variation $J$, those eigenfunctions $u$ simply loose contrast whereas the overall structure of the function remains the same. The magnitude of contrast loss is related to the eigenvalue $\lambda$. The eigenfunctions shape is determined by the chosen norm in the TV functional which can be adapted to the application of interest. In this way signals are not linearly smoothed to overcome scales but are step-by-step transformed to a composition of nonlinear eigenfunctions at coarser scales. A spectral response function can be used to examine which scales have a strong contribution to the original signal and to design filters of certain scales. Moreover, such spectral method can be combined with forward and inverse scale space approaches \cite{Burger2016}.\newpage

\textit{Main contribution.} In this paper we will extend the idea of (inverse) scale space methods known for nonlinear diffusion processes to segmentation and shape detection problems.
\begin{itemize}
\setlength{\itemsep}{3pt}
	\setlength{\itemindent}{-.2in}
	\item First goal: A novel inverse multiscale segmentation method based on Bregman iterations
	\item Second goal: An adaptive regularization parameter strategy for $\alpha$ (independent of $c_1$, $c_2$)
	\item Third goal: A spectral analysis for segmentation regarding shapes of eigenfunctions
\end{itemize}
A TV-based forward scale space for segmentation can easily be derived from the CV model via an increasing regularization parameter. We extend this framework by an inverse scale space for segmentation, still based on the CV model and therefore on a nonlinear diffusion process. For this purpose, we make use of Bregman iterations, among others well-known for improving total variation denoising results. The relation between the forward scale space and the inverse scale space is examined. Both iterative strategies are accomplished by a spectral transform and response function, which are used to easily examine scales and to filter certain scales. Since our method uses the total variation as a nonlinear diffusion process, we can make use of relatively easy and fast numerical and parallel implementation schemes developed in recents years.

\textit{Organization.} This work is organized as follows.
In Section \ref{sec:CVmodel} we start with a revision of the segmentation model by Chan and Vese including its convexification. Together with a revision of the iterative denoising strategy using Bregman iterations in Section \ref{ref:bregmandnoising}, the combination of those two concepts forms the first ingredient of our new inverse scale space method for multiscale segmentation introduced in Section \ref{sec:bregman-cv}. An interesting interpretation of the method as an adaptive regularization method is presented in Section \ref{sec:reg-param}. 
Section \ref{sec:spec-analysis} deals with nonlinear spectral methods and contains the second main ingredient of our new approach. In \ref{sec:genTV} we start with a brief summary of generalized nonlinear total variation functionals addressing different eigenshapes and continue in \ref{sec:spec-analysisdenoising} with recent works by Gilboa et al. solving related nonlinear eigenvalue problems in imaging. In Section \ref{sec:spec-analysissegm} we extend those ideas from nonlinear image denoising to image segmentation.
In Section \ref{sec:numrealization} we describe the numerical realization of our approach using primal-dual convex optimization methods. 
In Section \ref{sec:results} we first illustrate the strengths and limitations of our multiscale segmentation method by studies on synthetic datasets with a certain focus on eigenfunctions and shapes. Moreover, we underline the potential and wide applicability by three different biomedical imaging applications. Its reliable performance is demonstrated on real fluorescence microscopy images that contain Circulation Tumor Cells, with various shapes and sizes, among white blood cells and debris in \ref{sec:resultsCTC}. Besides, we present results on electron microscopy images suffering from inhomogeneous backgrounds in \ref{sec:resultsEM} and interesting results on network-like shapes representing vascular systems in \ref{sec:resultsnetwork}.
We end with a conclusion and an outlook to future possible perspectives in Section \ref{sec:conl}.
%
%
%
%
\section{Modeling Segmentation with Inverse Scale Spaces}%
In the following section we will shortly describe the model by Chan and Vese \cite{Chan2001} for segmentation and the adaption of the ROF model for denoising using Bregman distances \cite{Osher2005}. Afterwards we will introduce our novel Bregman-CV model for segmentation and show some advantages of our model. 
\subsection{Globally Convex Segmentation}\label{sec:CVmodel}
The idea of the CV model has originally been derived from the more general model for image segmentation introduced by Mumford and Shah (\cite{Mumford1989}). Here, one seeks for a solution of the variational energy
\begin{equation*}
J^{\text{MS}}(u,C) = \int_{\Omega} |f(x) - u(x)|^2 \mathrm{d}x +  \alpha \cdot Per(C) + \beta \cdot \int_{\Omega\setminus C} |\nabla u(x)| ^2 \mathrm{d}x \longrightarrow \min_{u,C}
\end{equation*}
where $u$ is a differentiable function that is allowed to be discontinuous on $C$. $C$ describes the union of the boundaries and thereby represents the contour defining the segmentation. Thus, $u$ is a smooth approximation of the original image $f$ and is composed of several regions $\Omega_i$. Within each region $\Omega_i$, $u$ is smooth. If we restrict this model so that $u$ is composed of only two regions $\Omega_1$ and $\Omega_2$ we derive with $\beta \rightarrow \infty$ the CV model for segmentation
\begin{equation*}
J^{\text{CV}}(c_1,c_2,C) = \int_{\Omega_1} (f(x) - c_1)^2 \mathrm{d}x + \int_{\Omega_2} (f(x) - c_2)^2 \mathrm{d}x + \alpha \cdot Per(C)\longrightarrow \min_{C,c_1,c_2}.
\end{equation*}
Here, $c_i$ is the intensity value of $u$ within the corresponding region $\Omega_i$. These values need to be determined together with contour the $C$. Note that in comparison to the original model we omit the area regularization term (cf. \eqref{eq:CVorig}).

The contour $C$ can be indirectly represented via a level-set function $\Phi$ (\cite{Osher1988}) with
\begin{equation*} 
C = \{x \in \Omega: \Phi(x) = 0\} 
\end{equation*} 
and $\Phi(x)$ being positive if and only if $x \in \Omega_1$. Together with the Heaviside function
\begin{equation*}
	H(\Phi(x)) = \begin{cases}1 &\text{ if }\Phi(x)\geq 0 \\ 0 &\text{ if }\Phi(x)< 0\end{cases}
\end{equation*}
and its regularized version $H_{\epsilon}$ this results in
\begin{align*}
J^{\text{CV2}}(c_1,c_2,\Phi) = \int_{\Omega} (f(x) - c_1)^2 H_{\epsilon}(\Phi(x)) \mathrm{d}x & + \int_{\Omega} (f(x) - c_2)^2 (1 - H_{\epsilon}(\Phi(x)))\mathrm{d}x \\&+ \alpha \cdot \int_{\Omega}|\nabla H_{\epsilon}(\Phi(x))|\mathrm{d}x~\longrightarrow~ \min_{\Phi,c_1,c_2}.
\end{align*}
The contour $C$ evolves during minimization until it reaches a minimum which, in the ideal case, describes the object boundaries. Besides the original minimization strategy by gradient descent, several minimization methods to solve the CV model have been developed, see for example \cite{He2007,Zehiry2007,Badshah2008,Bae2009}. One disadvantage of the model is its non-convexity which makes the solution depending on the used initialization. With a badly chosen initialization the minimization might get stuck in a local minimum that corresponds to a bad or meaningless segmentation.

For a better understanding of the relation between the nonlinear denoising model by Rudin, Osher and Fatemi (\cite{Rudin1992}) and the CV segmentation model we use the total variation defined as
\begin{equation}\label{eq:TV}
TV(u) := \sup\limits_{\substack{\varphi \in C_0^{\infty}(\Omega; \mathbb{R}^2)\\ ||\varphi||_{\infty} <1}} \int_{\Omega} u \nabla\cdot \varphi \;\mathrm{d}\mu \text {\ \  with \ \ } BV(\Omega) := \{ u \in L^{1}(\Omega)|TV(u) < \infty\}.
\end{equation}
The total variation of a characteristic $u(x) = \begin{cases} 1 &\text{\  if  \ } x \in \Omega_1 \cup C\\ 0  &\text{ if } x \in \Omega_2\end{cases}$ corresponds to the contour length $|C|$ which can be shown by the co-area-formula.
Therefore we can formulate the segmentation problem as
\begin{equation} \label{eq:CV3} J^{CV3}(c_1,c_2,u) = \int_{\Omega} u\left((f(x)-c_1)^2 - (f(x) - c_2)^2\right)\mathrm{d}x + \alpha \; TV(u)
\longrightarrow \min_{\substack{u\in BV(\Omega),c_1,c_2\\u(x)\in \{0,1\}}}.
\end{equation}
For fixed $c_1, c_2$ the solution of \eqref{eq:CV3} corresponds to the solution of an ROF problem with binary constraint (\cite{Burger2012}):
\begin{equation}\label{eq:CV-ROF}
 \min_{\substack{u\in BV(\Omega)\\u(x)\in \{0,1\}}}\frac{1}{2}|| u(x) - r(x)||_2^2 + \alpha  TV(u)
 \end{equation}
with $r(x) = (f(x)-c_2)^2 \!-\! (f(x) - c_1)^2 \!-\!\frac{1}{2}$.

The regularization parameter $\alpha$ in the segmentation model \eqref{eq:CV3} has the role of a scale parameter, meaning that $\alpha$ determines the scale of the objects that are segmented. The CV model describes a forward scale approach, thus a small parameter $\alpha$ corresponds to small scales that are segmented. An increased regularization parameter results in a solution where the smaller scales are not segmented but only larger ones. The meaning of scale is determined by the regularization functional, in this case the total variation. The total variation encodes a measure of the contour length as well as the height of piecewise constant areas. One disadvantage is that due to the 1-homogeneity of $TV$ our method cannot distinguish between height and contour length. Thus, a small object with a bright intensity can have the same scale as a large object with a less bright intensity. For more details see section \ref{sec:results}.
\
\paragraph{Convexification}
The CV segmentation model \eqref{eq:CV3} as well as the binary ROF model \eqref{eq:CV-ROF} are both not convex. Even for fixed values of $c_1$ and $c_2$ both models are non-convex due to the binary constraint on $u$. As mentioned before this might result in local instead of global minimum solutions. Approaches to overcome this difficulty and find global minima of the CV model are presented for example in\cite{Chan2006,Bresson2007,Goldstein2010,Brown2012}. In \cite{Chan2006} the authors showed that global minimizers of \eqref{eq:CV3} for any given fixed $c_1,c_2 \in \mathbb{R}$ can be found by solving
 \begin{equation}\label{eq:convCV}
J^{CV3}(c_1,c_2,v) = \int_\Omega v ((f(x)-c_1)^2 - (f(x)-c_2)^2) dx + \alpha~TV(v) \longrightarrow \min_{v \in BV(\Omega),~v(x) \in [ 0,1 ]}
\end{equation}
and defining $ u(x) := v^{\ast}(x) \geq \mu$ for a.e. $\mu \in [0,1]$.

Thus, the binary constraint can be relaxed and combined with a thresholding. Here, the variational model to solve is convex, though not strictly convex. One should bear in mind that the found solution is therefore not unique. Yet solutions of \eqref{eq:convCV} are close to binary even if the constraint is relaxed. Thus for most choices of $\mu$ we derive the same solution which means that the choice of $\mu$ has only a very limited impact on our method. Therefore we don't see any disadvantages when choosing the global but not unique minimum $u(x) := v^{\ast}(x) \geq 0.5$. A fully convex formulation (including the constants $c_1$ and $c_2$) of problem \eqref{eq:CV3} can be found in \cite{Brown2012}. This method is computationally less efficient and currently we don't think that, in our method, the advantages of the full convexity outweigh the increased computational time. 
\subsection{Inverse Scale Space for TV-Denoising}\label{ref:bregmandnoising}
Before introducing our new segmentation model in the following subsection we will first recall some properties of the well-known ROF model \cite{Rudin1992} and its extension by Bregman distances introduced in \cite{Osher2005}. To denoise an image corrupted by additive Gaussian noise, \cite{Rudin1992} proposed to solve the nonlinear variational problem  
\begin{equation}\label{eq:ROF}
\frac{1}{2}||u - f||_2^2 + \alpha \ TV(u) \longrightarrow \min_{u \in BV(\Omega)}
\end{equation}
referred to as the ROF model. 
Similar to the CV model this generates a forward scale space flow regarding the scale parameter $\alpha$. An increased parameter $\alpha$ leads again to a solution $u$ where fine scales are removed and vice versa. The total variation regularization functional is known to lead to a systematic contrast loss in the denoised image $u$ \cite{Meyer2001}. To tackle this problem Osher et al. proposed in \cite{Osher2005} an iterative contrast enhancement procedure based on Bregman distances. Instead of using the total variation regularization functional as before, information about the solution $u$ that we gained from a prior solution of problem \eqref{eq:ROF} is included. Therefore, problem \eqref{eq:ROF} is replaced by a sequence of variational problems 
\begin{equation}\label{eq:Bregman-ROF}
u_{k+1} = \argmin_{u\in BV(\Omega)} \frac{1}{2}||u - f||_2^2 + \alpha\ D^{p_k}_{TV}(u,u_k).
 \end{equation}
The regularization $D^{p_k}_{TV}(u,u_k):=TV(u)-TV(u_k)-\langle p_k,u-u_k\rangle$ is the Bregman distance of $u$ to the previous iterate $u_k$ with respect to the total variation. $p_k \in \partial TV(u_{k})$ is an element in the subdifferential of the total variation of the prior solution $u_{k}$. Although this subdifferential might be multivalued, the iterative regularization algorithm automatically selects a unique subgradient based on the optimality condition. For $k = 0$ we set $u_0 = p_0 = 0$.
The iterative strategy of this model is as follows: We start with a large parameter $\alpha$ that results in an oversmoothed solution $u_1$ that consists of only large scales. In every iteration step finer scales are added back to the solution. Thus, the scale parameter that determines the range of the scales present in $u$ is the iteration parameter $k$. In contrast to the forward approaches presented before, a small $k$ corresponds to coarse scales and a large $k$ to very fine scales. Therefore, the Bregman-ROF denoising approach is an inverse scale space approach. The authors showed that this strategy leads to enhanced contrast of the final solution $u_{k_{\text{max}}}$ compared to the solution of \eqref{eq:ROF}. Hence, solving \eqref{eq:Bregman-ROF} instead of \eqref{eq:ROF} with increasing $\alpha$ is not only an inverted way of detecting scales. This method rather allows for a detection of solutions which cannot be obtained by an adequate choice of $\alpha$ in the original ROF model.\\

\subsection{Bregman-CV Segmentation Model}\label{sec:bregman-cv}
In the following section we will introduce our new inverse scale space approach for segmentation. It is based on the similarity of the ROF functional and the CV functional shown in \eqref{eq:CV-ROF}. Similar to the Bregman-ROF denoising problem we replace the total variation regularization by an iterative regularization based on Bregman distances. Thus, the resulting novel segmentation model is given by
\begin{equation*}
u_{k+1} = \argmin_{\substack{u\in BV(\Omega)\\u(x)\in [0,1]}} \int_{\Omega} u\left((f-c_1)^2 - (f-c_2)^2\right)
+ \alpha \ D^{p_k}_{TV}(u,u_k)
\end{equation*}
By inserting the definition of the Bregman distance $D^{p_k}_{TV}(u,u_k):=TV(u)-TV(u_k)-\langle p_k,u-u_k\rangle$ and ignoring the parts independent of $u$ we derive the following model.
\begin{equation}\label{eq:Bregman-CV}
u_{k+1}  = \argmin_{u\in BV(\Omega)} \int_{\Omega} u\left((f-c_1)^2 - (f-c_2)^2\right) + \chi_{[0,1]}(u)
+ \alpha \ \left(TV(u)-<u,p_k>\right)
\end{equation}
with $p_k \in \partial TV(u_k)$, $p_0 = 0$ and $ \chi_{[0,1]}(u) = 0$ if $u(x) \in [0,1]$ and equal to infinity elsewhere. The range of the scales present in $u_{k+1}$ is again determined by the iteration index $k$. This model is an inverse scale space approach, thus a small $k$  corresponds to a large scale segmentation and vice versa. 

By definition of the Bregman distance it is $p_k \in \partial TV(u_k)$ where the subdifferential is multivalued. Therefore we need to determine a rule to choose $p_k$. One way is to derive an update strategy based on the optimality condition of \eqref{eq:Bregman-CV} (cf. \cite{Osher2005}):
\begin{equation*}
0 = \left((f-c_1)^2 - (f-c_2)^2\right) + q_{k+1} + \alpha p_{k+1} - \alpha p_k\text{\ \ \ (opt.cond.)}.
\end{equation*}
Here, $q_{k+1}$ is an element in the subdifferential of the characteristic function $\chi_{[0,1]}(u_{k+1})$. The subdifferential of a characteristic function is a normal cone and is in our case given by 
\begin{equation*}
q_k(x) \in \begin{cases} (-\infty,0] &\mbox{if} \quad u_k(x) = 0 \\ \{ 0 \} &\mbox{if} \quad 0 < u_k(x) < 1 \\ [0,\infty) &\mbox{if}\quad u_k(x) = 1\end{cases}.
\end{equation*}
Thus, we can choose $q_k = 0$ and neglect it from hereon. The update strategy for $p_{k+1}$ is then given by
\begin{equation}
p_{k+1} = p_k - \frac{1}{\alpha}\left((f-c_1)^2 - (f-c_2)^2\right) = -\frac{k+1}{\alpha}\left((f-c_1)^2 - (f-c_2)^2\right)\label{eq:update_p}.
\end{equation}
This update is independent of $u_k$, thus it is a pointwise constant update in every iteration.

\subsection{Interpretation as an Adaptive Regularization Approach}\label{sec:reg-param}
One important question is whether solutions of this model are in some sense improved compared to the solutions of the original CV model. We mentioned before that Bregman iterations lead to a contrast enhancement when applied to the ROF functional and that solving the CV model corresponds to solving a binary ROF. Yet, one should bear in mind that a contrast enhancement is meaningless in the case of a binary image since the contrast is already determined by the binary constraint. This is supported by the following observation: when inserting \eqref{eq:update_p} into \eqref{eq:Bregman-CV} we get
\begin{align*}
u_{k+1} & = \argmin_{\substack{u\in BV(\Omega)\\u(x)\in [0,1]}}\int_{\Omega} u\left((f-c_1)^2 - (f-c_2)^2\right) + \alpha \left(TV(u)-<u,p_k>\right)\\
& =  \argmin_{\substack{u\in BV(\Omega)\\u(x)\in [0,1]}} \int_{\Omega} u\left((f-c_1)^2 - (f-c_2)^2\right) + \frac{\alpha}{(k+1)} \ TV(u)
\end{align*}
With this, it is straightforward to see that all solutions $u \in BV(\Omega)$ derived by the Bregman-CV model can also be found by the original CV model. Nevertheless, there are advantages of using the iterative update strategy.
\begin{figure}[t]
  \centering 
  \subfigure[Linearly spaced $\alpha$'s from 50 to 1.]{\includegraphics[width=0.35\textwidth]{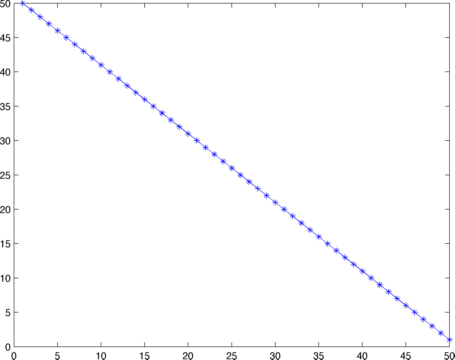}}\label{fig:reg_params_linear }\quad 
  \subfigure["$\alpha$'s" resulting from 50 Bregman steps with $\alpha = 50$.]{\includegraphics[width=0.35\textwidth]{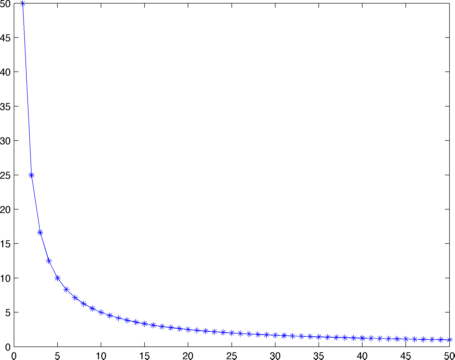}}\label{fig:reg_params_bregman}\quad 
  \caption{Comparison between linearly spaced regularization parameters and the automatically chosen parameters in the Bregman-CV model.}
\label{fig:reg_params} 
\end{figure}
In Figure \ref{fig:reg_params} b) the resulting $\tilde\alpha = \frac{\alpha}{k+1}$ for $\alpha = 50$ and 50 Bregman iterations are shown. It is obvious that in the first iterations the decrease of the regularization parameter is much larger compared to later iterations. This is reasonable since first the large scales are reconstructed and in later Bregman iterations the finer scales are incorporated in the result. Making large steps with $\alpha$ in large scales and becoming finer is therefore a reasonable strategy. A large decrease in the later iterations would probably miss some scales in between while in the first iterations too small steps are not reasonable. With this strategy the problem of automatically choosing the regularization parameter is less severe. By choosing a large $\alpha$ and performing multiple Bregman iterations, a broad spectrum of scales is detected. Yet one should bear in mind that a strategy to automatically detect important scales is needed for a fully automated framework. 
%
%
%
%
\section{A Spectral Method for Multiscale Segmentation}\label{sec:spec-analysis} 
The analysis of eigenvalues and spectral decomposition is a well-known theory in the field of linear signal and image processing (see e.g. \cite{MarpleJr1987} or \cite{Stoica2005} for a more recent overview). Since nonlinear regularization became popular in the last years, there is a growing interest in generalizing this theory to nonlinear operators. In \cite{Benning2012} Benning et al. examined singular values for nonlinear, convex regularization functionals and in  \cite{Gilboa2013,Gilboa2014} Gilboa transferred the idea of spectral decompositions to the nonlinear total variation functional and related operators. The general idea is to examine solutions of the nonlinear eigenvalue problem 
\begin{equation}\label{eq:eigenval}
\lambda u \in \partial J(u)
\end{equation}
where $\partial J$ denotes the subdifferential of a (one-homogeneous) convex functional usually representing regularization in inverse imaging problems. 
By transferring solutions of \eqref{eq:eigenval} to sparse peaks in a spectral domain, advanced filters enhancing or suppressing certain image components can easily be designed. This concept was first introduced for the total variation functional and later generalized to one-homogenous functions and analyzed in \cite{Burger2015}, \cite{Gilboa2015} and \cite{Burger2016}.

\begin{figure}[t]
  \centering 
  \subfigure[Solution of the ROF model for different $\alpha$.]{\includegraphics[height=0.22\textheight]{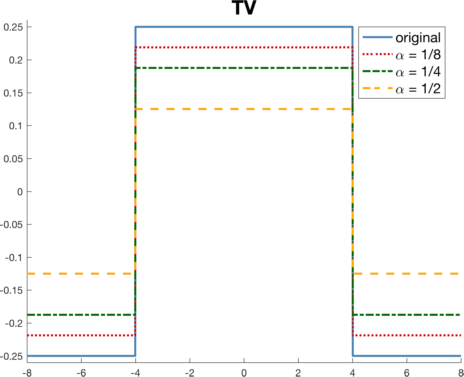}}\label{fig:TVblock }\quad 
  \subfigure[Solution of the inverse Bregman-ROF model for different Bregman iterations.]{\includegraphics[height=0.22\textheight]{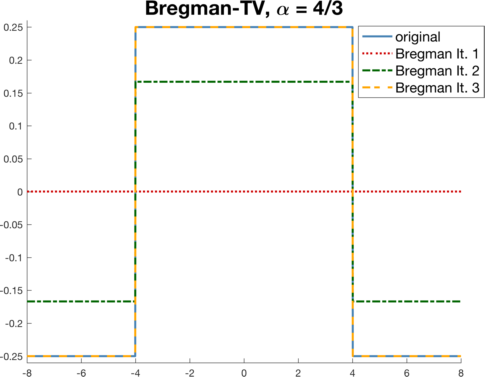}}\label{fig:BregmanTVblock}\quad 
  \caption{Solutions of the ROF and Bregman-ROF model in case the signal is a TV eigenfunction. Different regularization parameters $\alpha$ in the forward scale (left) and Bregman steps in the inverse scale space (right).}
\label{fig:1dEF} 
\end{figure}
In the case of a one dimensional signal the eigenfunction of the total variation corresponds to a single block with constant height, see Figure \ref{fig:1dEF}. When increasing the regularization parameter in the ROF model \eqref{eq:ROF}, i.e. when minimizing the total variation of this block signal, the block looses it height but the edges remain at the same position (cf. Figure \ref{fig:1dEF} (a)). Therefore all signals can be seen as the original signal multiplied by a scalar and are eigenfunction of the total variation. The same holds true for solutions of the Bregman-ROF model \eqref{eq:Bregman-ROF}. The only difference is that now eigenfunctions of the TV functional (blocks) are reshaped with increasing number of Bregman iterations instead of removed (see Figure \ref{fig:1dEF}(b)). This reflects the inverse scale approach of the Bregman-ROF model. Different to our novel approach, the Bregman updates in the Bregman-ROF model cannot be reformulated as an adaptive regularization parameter choice. As far as we know it is not entirely clear how the forward and inverse TV flow for denoising relate to each other.

\subsection{Generalized Definition of the Total Variation}\label{sec:genTV}
In the previous paragraph we mentioned that in the one dimensional case the eigenfunction of the total variation functional is a block with constant height. If we want to proceed to the two dimensional case  the eigenfunctions are not as clearly defined as in the 1D case. Different than before we now have the freedom to choose the norm body that is used within the infinity norm in the TV definition \eqref{eq:TV}. This choice determines the shape of the eigenfunctions. To reflect this dependence on the chosen norm in the TV definition we introduce a generalized version of the TV definition:
\begin{equation}\label{eq:TVgen}
TV_\gamma(u) :=  \sup\limits_{\substack{\varphi \in C_C^{1}(\Omega; \mathbb{R}^d)\\ \varphi(x) \cdot n <\gamma(n) \forall n\in\mathbb{R}^{d}}}\ -\int_{\Omega} u \nabla\cdot \varphi  \mathrm{d}x.
\end{equation}
Here, $\gamma : \mathbb{R}^d \rightarrow \mathbb{R} $ is a convex, positively 1-homogeneous function such that  $\gamma (x) > 0$ for $x \neq 0$.
If $u$ is a function in $\W^{1,1}(\Omega)$ the primal definition of equation \eqref{eq:TVgen} is given by
\begin{equation*}
TV_\gamma(u) :=  \int_{\Omega} \gamma(\nabla u) \mathrm{d}x.
\end{equation*}
In both definitions the choice of $\gamma$ determines the shape of the eigenfunctions. We refer to 
 \begin{equation*}
 \mathcal{F}_{\gamma} := \large\{ z \in \mathbb{R}^{d} : \gamma(z) \leq 1\large\}
 \end{equation*}
as the Frank diagram and the corresponding Wulff shape is defined as
\begin{align*}\mathcal{W}_{\gamma} :=& \large\{ z \in \mathbb{R}^{d} : z \cdot x \leq \gamma(x) \text{ for all } x\in \mathbb{R}^{d} \large\}\\
= &\large\{ z \in \mathbb{R}^{d} : \gamma^{\ast}(x) := \sup\limits_{x\in \mathbb{R}^{d}} \ \frac{z \cdot x}{ \gamma(x)} \leq 1 \large\}.\end{align*}
Note that from definition \eqref{eq:TVgen} together with the definition of the Wulff shape we can conclude that $\varphi(x) \in \mathcal{W}_{\gamma}$. Therefore, for every choice of $\gamma$, functions with the same shape as the Wulff shape are eigenfunctions of $TV_\gamma$ (see \cite{esedoglu2004} Theorem 4.1). The Frank diagram $\mathcal{F}_{\gamma}$ and the Wulff shape $\mathcal{W}_{\gamma}$ for the three most common choices of $\gamma$ are presented in Table \ref{tab:shapes}. For simplicity we will omit the subscript $\gamma$ in $TV_\gamma$ from hereon, but we will come back to it in the numerical results in section \ref{sec:results}.
\begin{table}
\centering
\begin{tabular}{|c|c|c|}
\hline
$\gamma(x)$ & $\mathcal{F}_{\gamma}$ & $\mathcal{W}_{\gamma} $\\[7pt]
\hline
$\gamma = \| \cdot \|_{1}$ & \includegraphics[height = 0.05\textheight]{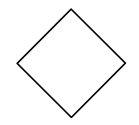} & \includegraphics[height = 0.05\textheight]{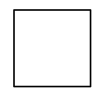}\\[6pt]
\hline
$\gamma = \| \cdot \|_{2}$& \includegraphics[height = 0.05\textheight]{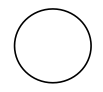} & \includegraphics[height = 0.05\textheight]{img/circle.png}\\[6pt]
\hline
 $\gamma= \| \cdot \|_{\infty}$ &\includegraphics[height = 0.05\textheight]{img/square.png} & \includegraphics[height = 0.05\textheight]{img/square_rot.png}\\
\hline
\end{tabular}
\caption{Examples for $\gamma$ and the corresponding Frank diagrams $\mathcal{F}_{\gamma}$ and Wulff shapes $\mathcal{W}_{\gamma}$. The Wulff shape corresponds to the shape of the eigenfunctions of the $TV_{\gamma}$ functional.}
\label{tab:shapes}
\end{table}

\subsection{Spectral Analysis for Nonlinear Functionals}\label{sec:spec-analysisdenoising}
In the following section we will first review the basic ideas of spectral TV analysis before transferring those ideas to segmentation in the next section. The idea to decompose an image based on the basic TV elements was presented by Gilboa in \cite{Gilboa2013,Gilboa2014}. These basic TV elements, called eigenfunctions of the total variation functional, are all functions $u \in BV(\Omega)$ that solve the nonlinear eigenvalue problem
\begin{equation}\label{eq:eigenvalProb}
\lambda u \in \partial TV(u).
\end{equation}
In the previous section we already presented some examples of TV eigenfunctions (see e.g. Figure \ref{fig:1dEF} and Table \ref{tab:shapes}). A  more general description of these eigenfunctions in case $\gamma$ is isotropic is given in \cite{Bellettini2002}. Bellettini et al. showed that all indicator functions $ \mathbbm{1}_C(x)$ of a convex and connected set $C$ with finite perimeter $Per(C)$ which admit
\begin{equation}\label{eq:l2eig}
\esssup_{p \in \partial C} \kappa (p) \leq \frac{Per(C)}{Area(C)}
\end{equation}
where $Area(C)$ denotes the area of $C$ and $\kappa$ the curvature of $\partial C \in C^{1,1}$ are solutions of \eqref{eq:eigenvalProb} and therefore eigenfunctions. Obtaining an analog condition for anisotropic, smooth and strictly convex choices of $\gamma$ is straightforward but more challenging in the case of non-smooth or non-strictly convex choices of $\gamma$. Candidates of these shapes are presented in \cite{bellettini2001}. In \cite{esedoglu2004} Esedoglu and Osher gave an example of a TV eigenfunction for $\gamma = \| \cdot \|_{1}$ that is not a Wulff shape (see section \ref{sec:resultsshapes} for more details).

In order to detect eigenshapes at different scales in a given signal $f$ a scale space approach is needed. One way to define a scale space based on the total variation is given by the total variation flow \cite{Andreu2001,Andreu2002,Bellettini2002,Burger2007,Steidl2004}. The TV flow arises when minimizing the total variation with steepest descent method and is defined as 
\begin{equation}\label{eq:tvflow}
\begin{aligned}
u_t(t,x) &= -p(t,x) \hspace{15pt} \text{for } p(t,x) \in \partial TV(u(t,x))\\
u(0,x) &= f(x)
\end{aligned}
\end{equation}
with Neumann boundary conditions. For $f(x) = \mathbbm{1}_C(x)$, with $C$ defined as above, the unique solution of \eqref{eq:tvflow} is $u(t,x) = (1 - \lambda_C t) ^{+}  \mathbbm{1}_C(x)$ with $\lambda_C = \frac{Per(C)}{Area(C)}$ (cf. \cite{Bellettini2002}). Hence, the time derivative $u_t(t,x)$ is given by the original signal multiplied with a scalar and $u$ is an eigenfunction. To obtain a suitable framework to decompose or filter images based on these eigenfunction Gilboa proposed to define a spectral framework that transforms eigenfunctions to peaks in the spectral domain. If $u(t,x)$ is a solution of \eqref{eq:tvflow} the TV spectral transform and spectral response can be defined as
\begin{equation}\label{eq:tvtransform}
\phi (t,x) = u_{tt}(t,x)\cdot t \quad\text{and}\quad S(t) = || \phi(t,x) ||_{L^{1}(\Omega)}.
\end{equation}
Note, that there are alternative definitions of the spectral response presented in \cite{Burger2015}. 

Another approach to construct an forward scale space is based on the variational ROF problem. Instead of solving \eqref{eq:tvflow} the ROF model
\begin{equation}\label{eq:ROFscalespace}
 \min_{u\in BV(\Omega)}\frac{1}{2}|| u(x) - f(x)||_2^2 + t \cdot  TV(u).
\end{equation}
is solved for different regularization parameters. Hence, $t$ is the (artificial) time variable that determines the scale comparable to $t$ in the TV-flow approach \eqref{eq:tvflow}. One drawback of this approach is that there is no clear rule for the choice of different $t$'s. The spectral transform and response function can be equivalently defined as in \eqref{eq:tvtransform}. 

A third, but significantly different, scale space approach is an inverse scale approach. The inverse scale space flow is defined as
\begin{equation}\label{eq:inverseflow}
\begin{aligned}
p_s(s,x) &= f(x) - u(s,x) \hspace{15pt} \text{for } p(s,x) \in \partial TV(u(s,x))\\
u(0,x) &= p(0,x) = 0.
\end{aligned}
\end{equation}
Thus, the flow is now defined on the dual variable $p \in \partial TV(u)$ and $s$ is the time variable determining the scale. Note that in \cite{Burger2005} it was shown that the iterative Bregman-L2-TV model \eqref{eq:Bregman-ROF} can be associated with a discretization of \eqref{eq:inverseflow} for $\frac{1}{\alpha} \rightarrow 0$.
As \eqref{eq:inverseflow} is an inverse approach the time variable $t$ can be associated with $\frac{1}{s}$. In this case the spectral transform and response functions are defined as
\begin{equation}\label{eq:tvtransform2}
\phi (s,x) = u_{s}(s,x) \quad\text{and}\quad S(s) = || \phi(s,x) ||_{L^{1}(\Omega)},
\end{equation}
where $u(s,x)$ is the solution of \eqref{eq:inverseflow}. Note, that for small $s$ $\phi(s,x)$ now measures changes in the coarse scales. See \cite{Burger2015} for more details.

With all three approaches we are able to transform a signal to the spectral domain and detect different scales based on TV eigenfunctions. If we assume that $\phi(t,x)$ is integrable over time, the original signal $f$ can be reconstruct via
\begin{equation*}
f(x) = \int_0^{\infty} \phi(t,x) dt + \bar{f}
\end{equation*}
where $\bar{f}$ is the average of $f$. Filters can be defined with $\phi_{H}(t,x) = H(t)\phi(t,x)$ via 
\begin{equation*}
f_{H}(x) = \int_0^{\infty} \phi_{H}(t,x) dt + H(\infty)\bar{f}.
\end{equation*}

\subsection{Spectral Response of Multiscale Segmentation}\label{sec:spec-analysissegm} 
In section \ref{sec:bregman-cv} we presented two variational models to detect segmentations of a given image $f$ at different scales. To decompose the segmentation into different scales and detect important scales or clusters of scales in the segmentation we want to transfer the idea of spectral analysis based on the total variation (cf. sec. \ref{sec:spec-analysisdenoising}) to segmentation. Therefore we need to find a suitable transformation of the segmentation $u$ to the spectral domain and vice versa. Note that our goal is not to reconstruct the original signal $f$ or filtered versions of $f$ but the reconstructed function should be a segmentation itself. To do so we make use of the idea that eigenfunctions of the TV functional should be transformed to peaks in the spectral domain. In the following we will derive this spectral transform function for a forward scale space and an inverse scale space approach. To represent the forward scale space we associate the regularization parameter $\alpha$ in the convex version of the original model by Chan and Vese \eqref{eq:convCV} with the artificial time variable $t$. That means, we solve
 \begin{equation}\label{eq:convCV-scale}
\int_\Omega u ((f(x)-c_1)^2 - (f(x)-c_2)^2) dx + t \cdot TV(u) \longrightarrow \min_{u \in BV(\Omega),~u(x) \in [ 0,1 ]}
\end{equation}
for different $t$.This is comparable to the variational approach in \eqref{eq:ROFscalespace}. So far, we could not find a forward scale space representation that can be associated with a flow on $u$ comparable to the TV-flow. One difficulty is that the optimality condition of this model, i.e. $0 = (f-c_1)^2 - (f-c_2)^2 + \alpha p$ with $p \in TV(u)$, has no direct dependence on $u$. The inverse scale space representation is based on the Bregman-CV model we introduced in \eqref{eq:Bregman-CV}. The optimality condition in each step of the iterative Bregman-CV strategy is given as
 \begin{equation*}
 0 = \left((f-c_1)^2 - (f-c_2)^2\right) + \alpha \left( p_k - p_{k-1}\right) \text{ with } p_k \in TV(u_k) \ \forall \ k \\.
 \end{equation*}
The resulting equation
\begin{equation*}
\frac{p_k - p_{k-1}}{\frac{1}{\alpha}} = (f-c_2)^2 - (f-c_1)^2
 \end{equation*}
can be interpreted as a discretization with stepsize $\frac{1}{\alpha}$ of
\begin{equation}\label{eq:invsegslow}
\begin{aligned}
p_s(s,x) &=  (f(x)-c_2)^2 - (f(x)-c_1)^2\hspace{15pt} \text{with } p(s,x) \in \partial TV(u(s,x))\\
p(0,x) &= 0.
\end{aligned}
\end{equation}
Again, $s$ is the time variable that is inverse to the time variable $t$ in \eqref{eq:convCV-scale}. We refer to this flow as the inverse scale space segmentation flow. Note that within this flow description there is no direct dependence on $u$ but $u$ is only indirectly given by $p \in \partial TV(u)$. Yet, when looking for solutions of this flow at different times $s$, we solve the Bregman-CV model for multiple Bregman iterations and there the corresponding $u$ is available. 

\begin{figure}[t]
  \centering 
  \subfigure[Evolution of $u$ in \eqref{eq:convCV-scale}.]{\includegraphics[width=0.4\textwidth]{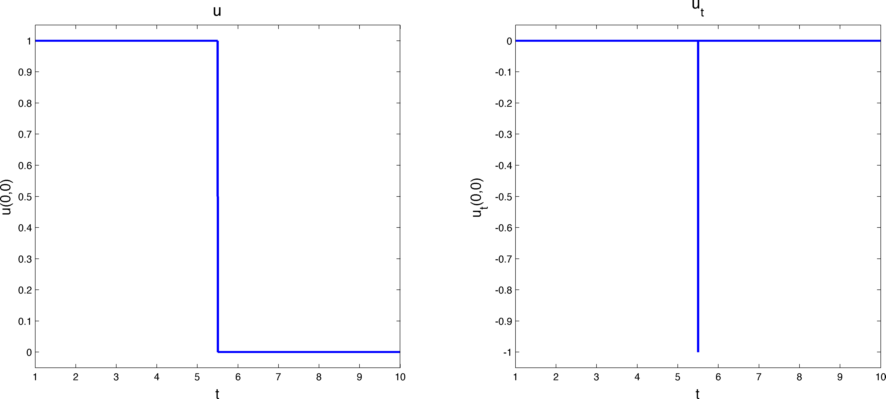}}\quad \quad
  \subfigure[Evolution of $u$ in \eqref{eq:invsegslow}.]{\includegraphics[width=0.4\textwidth]{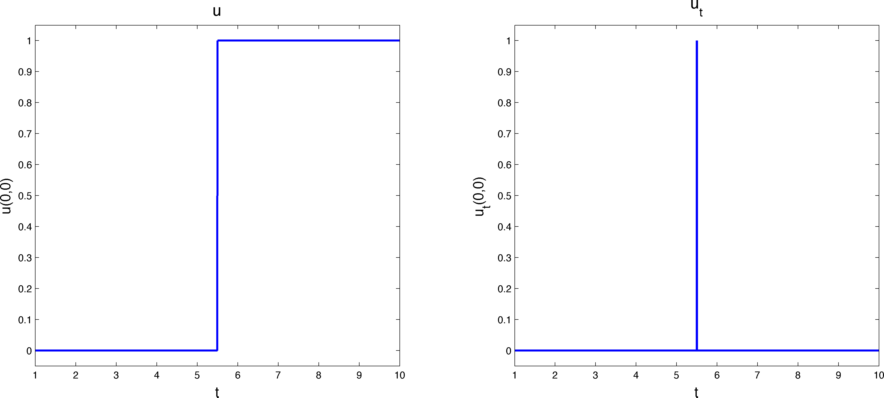}}
  \caption{Illustration of the evolution of $u\in\{0,1\}$ over time at a fixed point when $f$ is a TV eigenfunction shown for the forward (a) and inverse scale space approach (b). Note, due to the binary constraint on u the evolution of $u$ is already a step function, i.e. $u(x)$ is either part of the segmentation (u(x)=1) or not (u(x)=0).}
\label{fig:peaks} 
\end{figure}

To transform the eigenfunctions to peaks in the spectral domain we define the spectral transform function as follows:
\begin{equation}\label{eq:spectransseg}
\phi (t,x) =\begin{cases} -u_{t}(t,x)& \text{ (forward case)}\\
\ \ \ u_t(t,x)& \text{ (inverse case)}\end{cases}\\[3pt]
\end{equation}
This definition is motivated by Figure \ref{fig:peaks} where the evolution of TV eigenfunctions over time in \eqref{eq:convCV-scale} and \eqref{eq:invsegslow} is illustrated. We can see, that the first derivative (in a distributional sense) is a peak at that time point where the eigenfunction vanishes or appears respectively. For the spectral response function we use the definition
\begin{equation}
	S(t) = || \phi(t,x) ||_{L^{1}(\Omega)}
	\label{eq:spectralResponse}
\end{equation}
introduced by Gilboa (\cite{Gilboa2013,Gilboa2014}). Here, the influence of certain scales to the segmentation is encoded. Using this function $\phi$ and assuming integrability over time we can get the following backtransformation:
\begin{equation*}
f^{\text{seg}}(x) = \int_0^{\infty} \phi(t,x) dt,
\end{equation*}
where $f^{\text{seg}}(x):= (f - c_2)^2 - (f-c_1)^2 + \frac{1}{2} < \frac{1}{2}$ is the results of the simple clustering problem
$\min_{u(x)\in \{0,1\}} \int_{\Omega} u\left((f(x)-c_1)^2 - (f(x) - c_2)^2\right)\mathrm{d}x$. By choosing a function $H$ with $H(t) \in \{ 0,1 \}$ certain scales of interest can be filtered via
\begin{equation}\label{eq:filterseg}
f^{\text{seg}}_{H}(x) = \int_0^{\infty} \phi_{H}(t,x) dt \ \ \text{ with } \ \ \phi_{H}(t,x) = H(t)\phi(t,x).
\end{equation}
With this framework we can easily segment an image and simultaneously detect important scales in this segmentation instead of using the spectral TV analysis as in Section \ref{sec:spec-analysisdenoising} and segmenting separately. Moreover, using the filtering approach we can easily construct segmentations of only certain scales by filtering (cf. \eqref{eq:filterseg}). Examples are shown in section \ref{sec:results}.%
%
%
%
\section{Numerical Methods}\label{sec:numrealization}
Our novel approach introduced in the two previous sections consists of two parts. First we have to solve the Bregman-CV model  \eqref{eq:Bregman-CV} introduced in section \ref{sec:bregman-cv}. Afterwards we can analyze the detected scales using the spectral response function \eqref{eq:spectralResponse} introduced in section \ref{sec:spec-analysissegm}. Therefore an efficient solution of the Bregman-CV model is required. In the following section we will give a very brief introduction into primal-dual optimization schemes and then show how this can be used to solve \eqref{eq:Bregman-CV}. We close the section by a speed comparison between a Matlab and a parallelized C/Mex implementation of our code. 

\subsection{Primal-Dual Optimization Methods}
To solve nonlinear problems of the form
\begin{equation}\label{eq:primalprob}
 \min_{x \in X} \ F(Kx) + G(x)
 \end{equation}
with $F$ and $G$ being proper, convex, lower-semicontinuous functions, primal-dual optimization methods became very popular in the last years. Instead of minimizing the primal problem  \eqref{eq:primalprob} they make use of the primal-dual formulation of this problem given by
\begin{equation}\label{eq:primaldualprob}
\min_{x\in X} \max_{y\in Y}\  \langle Kx,y\rangle + G(x) - F^{\ast}(y)
\end{equation}
with $G: X \rightarrow [0, +\infty]$ and $F^{\ast}:Y \rightarrow [0, +\infty]$ being the convex-conjugate of $F$. By updating in every iteration step a primal and a dual variable these methods are able to avoid some difficulties that arise when working on the purely primal or dual problem. One example is the minimization of variational methods with TV regularization. If the gradient is zero, the TV functional is not differentiable which leads to problems in purely primal minimization schemes like gradient descent. Some examples for primal-dual minimization algorithms are the PDHG algorithm \cite{Zhu2008}, a generalization of PDHG by Esser et al. \cite{Esser2010}, the Split Bregman algorithm by Goldstein and Osher \cite{Goldstein2009}, Bregman iterative algorithms \cite{Yin2008}, the second-order CGM algorithm \cite{Chan1999} and inexact Uzawa methods \cite{Zhang2011}.

In \cite{Chambolle2011} Chambolle and Pock proposed an algorithm which can be seen as a generalization of PDHG as well. This algorithm was originally proposed by Pock et al. in \cite{Pock2009} to minimize a convex relaxation of the Mumford-Shah functional. The efficient first-order primal-dual algorithm to minimize general problems of the form \eqref{eq:primaldualprob} is presented in Algorithm \ref{alg:cp}. 
$ (I + \sigma \partial F^{\ast})^{-1}$ and $(I + \tau \partial G)^{-1}$ are the resolvent operators of $F^{\ast}$ and $G$ respectively which are defined through 
\eq{y = (I + \tau \partial G)^{-1}(x) = \argmin_{y} \left\{ \frac{\| y-x\|^2}{2\tau} + G(y) \right\}.}
\begin{algorithm}
	\caption{First-order primal-dual algorithm by Chambolle and Pock (\cite{Chambolle2011})}
	\label{alg:cp}
	{\fontsize{10}{10}\selectfont 
	\begin{algorithmic}
		\State \textbf{Parameters:} $\tau, \sigma > 0$, $\theta \in [0,1]$
		\State \textbf{Initialization:} $n=0,\ x^0 \in X,\ y^0 \in Y, \bar{x}^0 = x^0$\\
		\State \textbf{Iteration:}\\
                \State  \textbf{for } $(n\geq 0)$\ \textbf{ do }
		\begin{enumerate}\itemsep5pt
			\item $y^{n+1} = (I + \sigma \partial F^{\ast})^{-1}(y^n + \sigma K \bar{x}^n)$.
			\item $x^{n+1} = (I + \tau \partial G)^{-1}(x^n - \tau K^{\ast} y^{n+1})$.
			\item $\bar{x}^{n+1} = x^{n+1} + \theta (x^{n+1} - x^{n})$.
			\item Set $n=n+1$.
		\end{enumerate}\\
               
		\State \textbf{end for}\\
		\State \Return $x^n$
	\end{algorithmic}
	}
\end{algorithm}

\subsection{Numerical Realization Bregman-CV}
In order to use Algorithm \ref{alg:cp} to solve the constraint problem \eqref{eq:Bregman-CV} we reformulate the problem into
\begin{equation}\label{eq:primBregCV}
 \int_{\Omega} u\left((f-c_1)^2 - (f-c_2)^2\right)
+ \text{id}_{[0,1]}(u) + \alpha \ \left(TV(u)-<u,p_k>\right) \longrightarrow \min_{u\in BV(\Omega)}
\end{equation}
with $p_k \in \partial TV(u_k)$ and $p_0 = 0$. $\text{id}_{[0,1]}(u)$ is the indicator function of the interval $[0,1]$ defined as 0 if $u \in [0,1]$ and $\infty$ otherwise. To derive an minimization strategy based on Algorithm \ref{alg:cp} we set $x = u$, $K(u) = \grad u$ and
\begin{equation*}
F(u) = \| u\|_{1} \text{ and } G(u) = \text{id}_{[0,1]}(u) +  \int_{\Omega} u\left( (f - c_1)^2 - (f - c_2)^2 - \alpha p_k\right).
\end{equation*}
The convex-conjugate of $F(u) = \| u\|_{1}$ is given by $F^{\ast}(p) = \delta_{P}(p)$ with $P = \left\{p: \| p\|_{\infty} \leq 1\right\}$ and
\begin{equation}
	\delta_{P}(p) = \begin{cases} 0 &\mbox{if } p \in P\\ \infty &\mbox{if } p\notin P\end{cases} .
\end{equation}
Together with \eqref{eq:primaldualprob} we derive the primal-dual variant of \eqref{eq:primBregCV}:
\eqn{\label{eq:dualBregCV}\langle \grad u,p\rangle + \text{id}_{[0,1]}(u) +  \int_{\Omega} u\left[ (f - c_1)^2 - (f - c_2)^2 - \alpha p_k\right] - \alpha \delta_{P}(p) \longrightarrow \min_{u}\max_{p}.}
The resolvent operators for $G$ and $F^{\ast}$ are defined through
\begin{align}
u = (I + \tau \partial G)^{-1}(\tilde{u}) &= \text{Proj}_{[0,1]}\left[\tilde{u} - \tau\left((f - c_1)^2 - (f - c_2)^2 - \alpha p_k\right)\right]\notag\\
& = \max\left(0,\min\left(1, \tilde{u} - \tau\left((f - c_1)^2 - (f - c_2)^2 - \alpha p_k\right)\right)\right)\label{eq:proxprimal}
\end{align}
and
\begin{equation}
p = (I + \sigma \partial F^{\ast})^{-1}(\tilde{p}) = \text{Proj}_{\left\{\left\{p: \| p\|_{\infty} \leq 1\right\}\right\}}\left(\tilde{p_{ij}}\right).\label{eq:proxdual}
\end{equation}
Note that the $L^{\infty}$ norm $\| p\|_{\infty}$ is in the discrete setting defined as $\| p\|_{\infty} = \max_{i,j} \large\{\gamma^{\ast}(p_{i,j})\large\}$ where the choice of $\gamma$ determines the shape of the eigenfunctions of the TV functional. For $\gamma = \| \cdot \|_{\ell^{p}}$ with $p = 1$ the unit ball defined by $$\left\{ (x,y) \in \Omega | \max\{|x|,|y|\} \leq 1\right\}$$ is an TV eigenfunction, for $p = 2$ the unit ball defined by $$\left\{ (x,y) \in \Omega | \sqrt{|x|^2+|y|^2} \leq 1\right\}$$ and for $p=\infty$ the unit ball defined by $$\left\{ (x,y) \in \Omega | (|x|+|y|) \leq 1\right\}.$$ However these are not the only eigenfunctions.
With \eqref{eq:proxprimal} and \eqref{eq:proxdual}, we derive the primal-dual algorithm presented in Algorithm \ref{alg:cpforcv} to minimize  \eqref{eq:Bregman-CV}. Note that we are not updating the constants $c_1$ and $c_2$, but start with a good estimate and leave it fixed. To a certain extend the varying regularization parameter can compensate for an error in those constants. Some examples are presented in Section \ref{sec:results}.
\begin{algorithm}
	\caption{First-order primal-dual algorithm to solve \eqref{eq:Bregman-CV}.}
	\label{alg:cpforcv}
	{\fontsize{10}{10}\selectfont 
	\begin{algorithmic}%
		\State \textbf{Parameters:} data $f$, reg. param. $\alpha \geq 0$, $\tau, \sigma > 0$, $\theta \in [0,1]$, 			   $maxIts \in \mathbb{N}, \ maxBregIts \in \mathbb{N}$
		\State \textbf{Initialization:} $l=1,\ u^0_1=0, \ p_0:=0, \bar{u}^0 = u^0$\\
		\State \textbf{Iteration:}\\
		\State \textbf{while } \big($k<maxBregIts$\big) \textbf{ do }
                \begin{enumerate}\itemsep5pt
			\item Set n=0.
		\end{enumerate}
		\begin{quote}
		 \State  \textbf{while } \big($n<maxIts$\big) \textbf{ do }\\
		\begin{enumerate}\renewcommand{\labelenumi}{\alph{enumi})}\itemsep5pt
			\item $p^{n+1} = \text{Proj}_{\left\{\left\{p: \| p\|_{\infty} \leq 1\right\}\right\}}\left( p^n + \sigma \grad \bar{u}^n \right)$.
			\item $u^{n+1} = \max\left(0,\min\left(1, u^n + \tau \div p^{n+1} - \frac{\tau}{\alpha}\left((f - c_1)^2 - (f - c_2)^2 - \alpha p_k\right)\right)\right)$.
			\item $\bar{u}^{n+1} = u^{n+1} + \theta (u^{n+1} - u^{n})$.
			\item Set $n=n+1$.
		\end{enumerate}\\
		\State \textbf{end while}\\
		\end{quote}
		\begin{enumerate}\itemsep5pt
		\setcounter{enumi}{1}
		\item Update $p_{k+1} = p_k + \frac{1}{\alpha}\left((f - c_2)^2 - (f - c_1)^2\right)$.
		\item Set $u^0_{k+1}=u^n_{k}$.
		\item Set $k=k+1$.
		\end{enumerate}\\
		\State \textbf{end while}\\
		\State \Return $u^0_{k}$
	\end{algorithmic}
	}
\end{algorithm}
\begin{figure}[t]
  \centering 
  \includegraphics[width=0.6\textwidth]{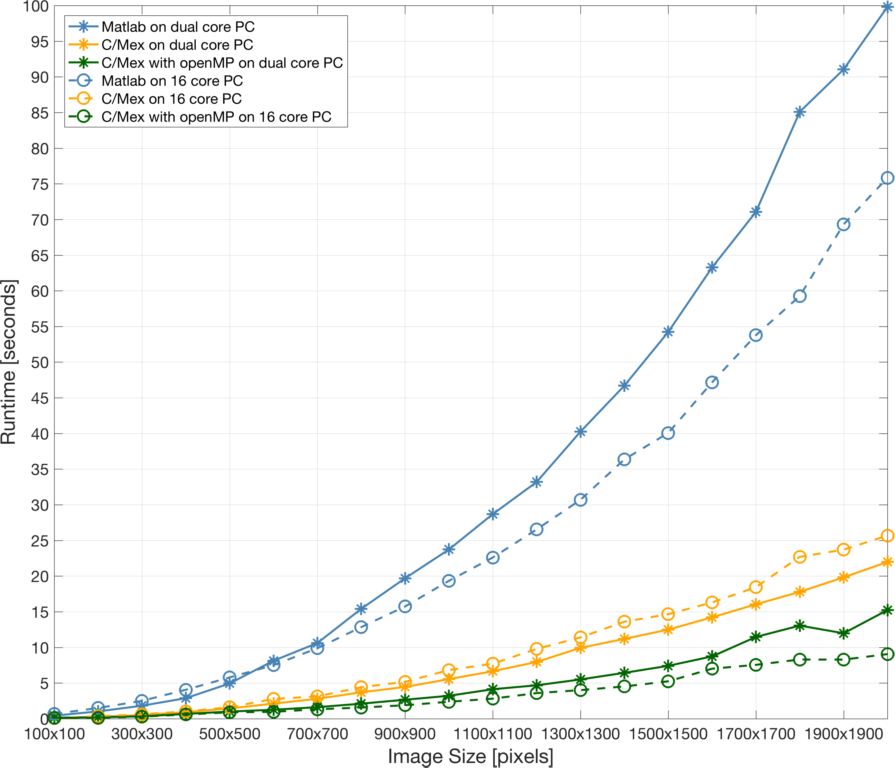}
  \caption{\textit{Comparison of runtimes for Matlab and C/Mex code.} The code was tested on 2 different computers: a Mac with an Intel Core i5 processor with two cores of 2.6 GHz and 8GB of memory (solid graphs) and a Linux computer with an Intel Xeon processor with 16 cores of  2.6GHz and 32 GB of memory (dashed graphs). On both computers we have tested a full Matlab implementation (blue), a C/MEX implementation without parallelization (yellow) and a C/MEX implementation with parallelization using openMP (green).}
\label{fig:speed} 
\end{figure}

We implemented two versions of the code, the first one is written in MATLAB and the second one is written in C and called via MATLAB using the MEX interface. For speed comparisons we've tested both code versions (always with 1 Bregman iteration) on two different machines. First, a MacBook Pro with an Intel Core i5 processor with two cores of 2.6 GHz and 8GB of memory (solid graphs in Figure \ref{fig:speed}) and the second machine was a Dell Linux workstation with in Intel Xeon processor with 16 cores of 2.6GHz and 32GB of memory (dashed graphs in Figure \ref{fig:speed}). Moreover, we ran the C/Mex code with and without parallelization using openMP. The results of our speed test can be seen in Figure \ref{fig:speed}. We see that especially for large images the computational time is significantly decreased when using the C implementation (yellow and green graphs) compared to the pure MATLAB code (blue graphs). All C versions behave very similar for small images but deviate for larger images. Without parallelization the C code is slightly faster on the Mac (solid yellow graph). When activating openMP and parallelizing parts of the code the computational time could be further improved for both computers. The difference is much more significant when using the computer with more cores. For the Mac we see a speed-up of about 30 percent (solid green graph) for large images while on the Linux computer (dashed green graph) we could achieve a speed-up of about 65 percent compared to the C implementation without parallelization.
%
%
%
%
\section{Experimental Results}\label{sec:results}
In this section we illustrate the main properties, advantages and limitations of our novel multiscale segmentation framework via several synthetic experiments as well as real experiments from biological cell imaging and retina imaging. The main framework, we apply and study throughout this section, consists of Algorithm \ref{alg:cpforcv}, performing the inverse scale space method Bregman-CV in \eqref{eq:Bregman-CV}, together with an evaluation of the spectral response $S$ in \eqref{eq:spectralResponse}.

In the first set of segmentation experiments we focus on simple objects (e.g. discs) and investigate how underlying object intensities and sizes can automatically be detected as multiple scales. In the second part we concentrate on the stability of our method under uncertainty like noise. In the third part of this section we generalize our study to other simple shapes (e.g. blocks or diamonds), explore the connection to underlying eigenfunctions of the nonlinear $TV_{\gamma}$ functional and the potential for shape detection respectively shape optimization. Finally, we apply our framework to real images obtained from the CellSearch system and other automatic scanning microscopes used for the enumeration of Circulating Tumor Cells \cite{allard2004} offering new automatic procedures for tumor cell identification, quantification and classification. Moreover, we show potential applications of our algorithm in the field of electron microscopy of biological samples and retina imaging.
\subsection{Multiscale Segmentation: Size versus Intensity}\label{sec:resultsscale}
In this first subsection we start with the investigation of basic segmentation tasks. We try to find round objects (discs) with varying scale (e.g. differences in size or intensity) in front of homogeneous backgrounds to verify the basic properties of our model. Since discs are the unit balls of the $\ell^2$ norm we chose $\gamma = \gamma^{\ast} = \| \cdot \|_2$ in the TV norm \eqref{eq:TVgen}.
\begin{figure}[htb]
  \centering 
    \subfigure[1st Bregman iteration.]{\includegraphics[width=0.18\textwidth]{results_new/balls_size_nonoise/bregman_cv_1.png}}\quad 
  \subfigure[9th Bregman iteration.]{\includegraphics[width=0.18\textwidth]{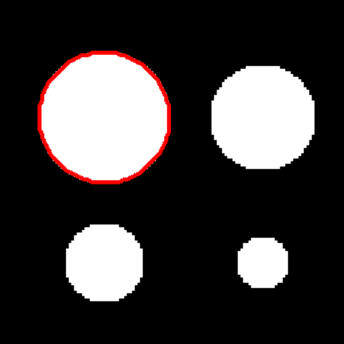}}\quad 
  \subfigure[12th Bregman iteration.]{\includegraphics[width=0.18\textwidth]{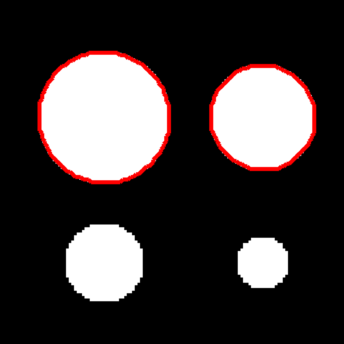}}\quad 
    \subfigure[15th Bregman iteration.]{\includegraphics[width=0.18\textwidth]{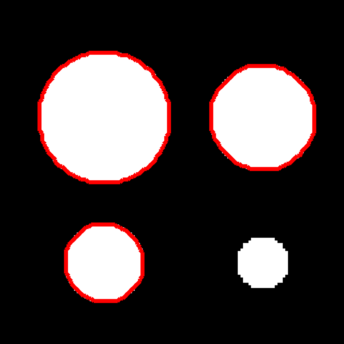}}\quad 
      \subfigure[22nd Bregman iteration.]{\includegraphics[width=0.18\textwidth]{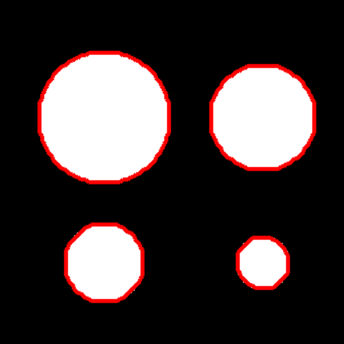}}\\
        \subfigure[Spectral response $S$]{\includegraphics[height=0.16\textheight]{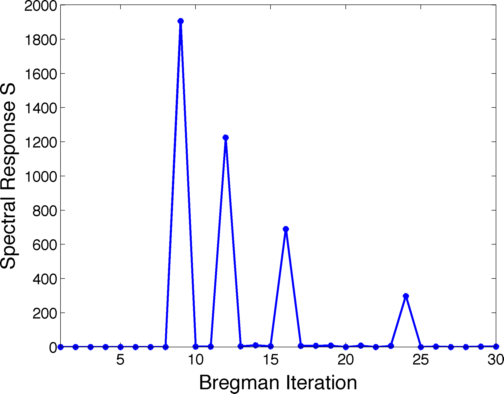}}\quad
         \subfigure[Color-coded segmentation of multiple scales.]{\includegraphics[height=0.16\textheight]{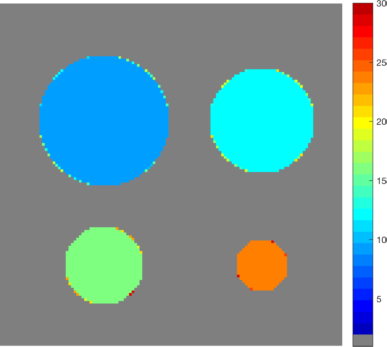}}
  \caption{\textit{Detection of size scales}. Automatic segmentation of discs of different sizes under fixed intensity (height). Segmentation is visualized via red contour. (a)-(e) Segmentation after different number of Bregman Iterations. (f) Associated spectral response function $S$. The number of Bregman iterations is plotted on the x-axis. Every peak in $S$ corresponds to one object that appears in the segmentation (compare to (a) - (e)). The decreasing size of the peaks reflects the decreasing disc size. }
\label{fig:ballssize} 
\end{figure}

In the first example in Figure \ref{fig:ballssize}(a) we can automatically identify four discs with different diameters (sizes) via the Bregman-CV inverse scale space algorithm \ref{alg:cpforcv}. As visualized in Figure \ref{fig:ballssize}(b)-(e) the method automatically picks out the discs one by one while iterating (w.r.t $t \approx \alpha_k$) through the scale space. At the end, the full foreground-background segmentation is available. In addition to this, the spectral response \eqref{eq:spectralResponse} of our iterates regarding $t$, nicely reflects the sparse occurrence of new contours, see Figure \ref{fig:ballssize}(f). Note that in this example the intensity of peaks decreases due to the different size of each disc. The height of $S$ at each time point reflects the number of pixels that occurred (or vanished) in the segmentation in this step. In Figure \ref{fig:ballssize}(g) we plotted the function $\int_{t}\Phi(x,t) \cdot t \mathrm{d} t$ in color coding. To obtain an easily interpretable visualization we use a gray background. Since the function $\Phi(x,t)$ for given $x$ is equal to 1 exactly at the time point $t_x$ where $x$ occurs in the segmentation and 0 elsewhere, $\int_{t}\Phi(x,t) \cdot t \mathrm{d} t$ corresponds for every point $x$ to $t_x$. We can see that apart from the four time steps where the four discs appear nearly nothing is appearing in between, reflecting the sparse behavior of our methods. Only very few pixels at the object boundaries occur in later iterations which is caused by discretization inaccuracies. In this way the multiscale approach of our method can be seen since each color reflects one segmented scale.\\
\begin{figure}[htb]
  \centering 
\subfigure[1st Bregman iteration.]{\includegraphics[width=0.18\textwidth]{results_new/balls_int_nonoise/bregman_cv_1.png}}\quad 
  \subfigure[3rd Bregman iteration.]{\includegraphics[width=0.18\textwidth]{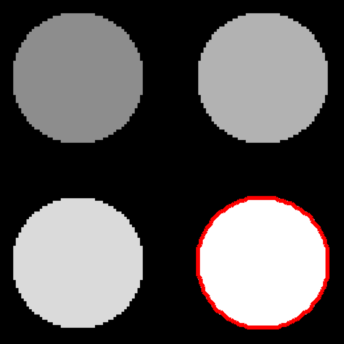}}\quad 
  \subfigure[5th Bregman iteration.]{\includegraphics[width=0.18\textwidth]{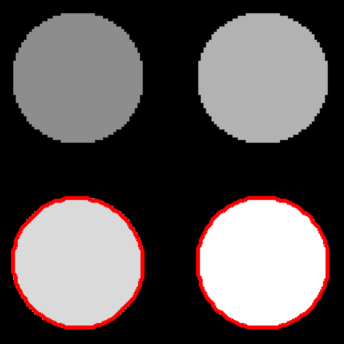}}\quad 
    \subfigure[10th Bregman iteration.]{\includegraphics[width=0.18\textwidth]{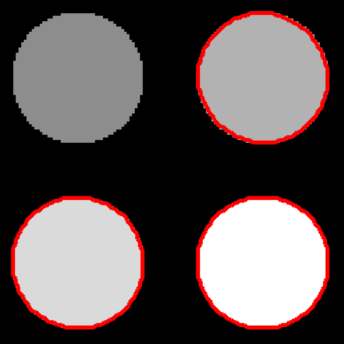}}\quad 
      \subfigure[29th Bregman iteration.]{\includegraphics[width=0.18\textwidth]{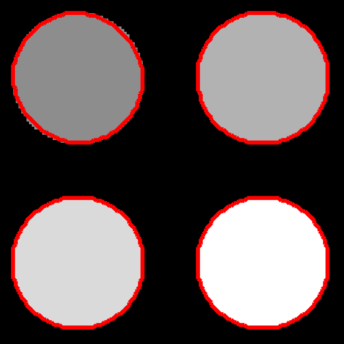}}\\
        \subfigure[Spectral response $S$.]{\includegraphics[height=0.16\textheight]{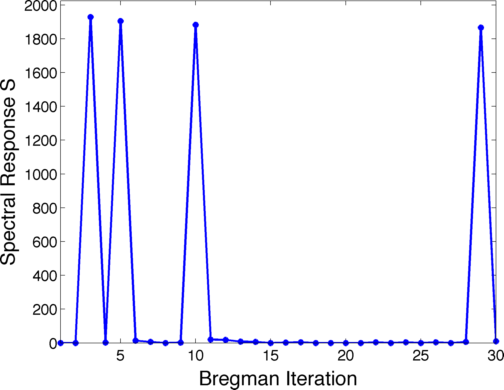}}\quad
         \subfigure[Color-coded segmentation of multiple scales.]{\includegraphics[height=0.16\textheight]{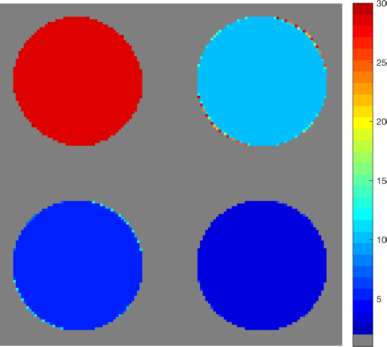}}

  \caption{\textit{Detection of intensity scales}. Automatic segmentation of discs of different intensities (heights) under fixed size. Segmentation is visualized via red contour. (a)-(e) Segmentation after different number of Bregman Iterations. (f) Associated spectral response function $S$. The number of Bregman iterations is plotted on the x-axis. Every peak in $S$ corresponds to one object that appears in the segmentation (compare to (a) - (e)).}
\label{fig:ballsint} 
\end{figure}

The results of the same experiment but now with intensity scales are shown in the next Figure.  The four discs in Figure \ref{fig:ballsint}(a) all have
the same size but vary in their intensity from $0.55$ to $1$. Note that in some sense this contradicts the assumption of the CV model that the data roughly consists of two different intensity levels $c_1$ and $c_2$. Nevertheless our multiscale approach is able to detect all scales reliably. Like in all our examples we estimated the constants beforehand as $c_1 = \text{mean}{\{f(x)| f(x) < \frac{1}{2}\cdot \max{f(x)}\}}$ and $c_2 = \text{mean}{\{f(x)| f(x) \geq \frac{1}{2}\cdot \max{f(x)}\}} $ and leave them fixed throughout the whole iteration process. Figure \ref{fig:ballsint}(b)-(e) visualizes the result of the four Bregman iterations in which the different discs appear. The color-coded result of multiples scales ($\int_{t}\Phi(x,t) \cdot t \mathrm{d} t$) is shown in Figure \ref{fig:ballsint}(g). Again, only a few pixels occur at time steps apart from the four significant ones. It is remarkable that the first three scales appear shortly after each other while the last scale appears strikingly later. This can be explained by the way the regularization parameter is adapted automatically through the Bregman distance (cf. Figure \ref{fig:reg_params} ). In later iterations the change in $\alpha_k$ becomes smaller so that the time until a certain object appears is larger compared to the previous iterate. Although this might lead to late appearance of smaller scales we can thereby guarantee to not miss scale differences in the finer scales. The late appearance of the last disc is also reflected in the spectral response function $S$ in Figure \ref{fig:ballsint}(f). Again, the response function has a sparse representation and the last peak occurs delayed. Apart from small variations caused by discretization artifacts the peaks now have all the same height reflecting their similar size.

\paragraph{Limitations:} In the previous two examples we showed how we can reliably detect different size scales and different intensity scales with our method. One limitation of our approach is the simultaneous detection of size and intensity scales in an image. 
\begin{figure}[htb]
  \centering 
  \subfigure[Data $f$ with intensity of 0.68 (large disc).]{\includegraphics[width=0.18\textwidth]{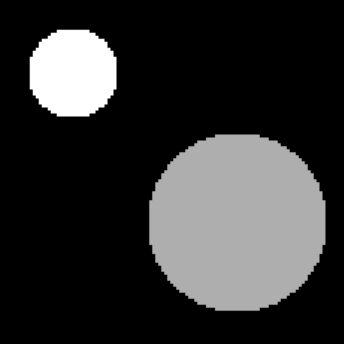}}\quad \quad
  \subfigure[Data $f$ with intensity of 0.69 (large disc).]{\includegraphics[width=0.18\textwidth]{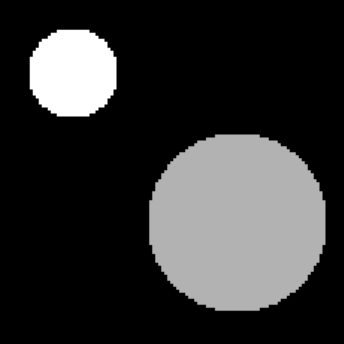}}\quad \quad
  \subfigure[Data $f$ with intensity of 0.70 (large disc).]{\includegraphics[width=0.18\textwidth]{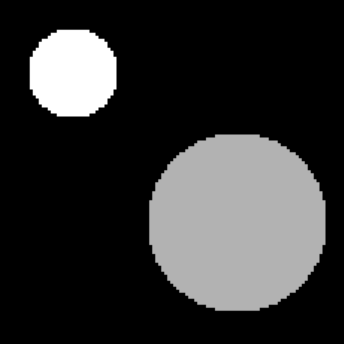}}\\
   \subfigure[Spectral response $S$ for (a).]{\includegraphics[width=0.2\textwidth]{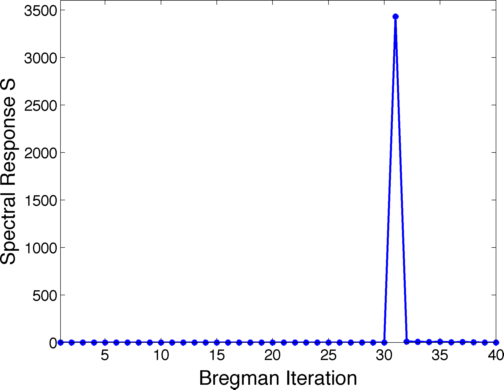}}\quad 
      \subfigure[Spectral response $S$ for (b).]{\includegraphics[width=0.2\textwidth]{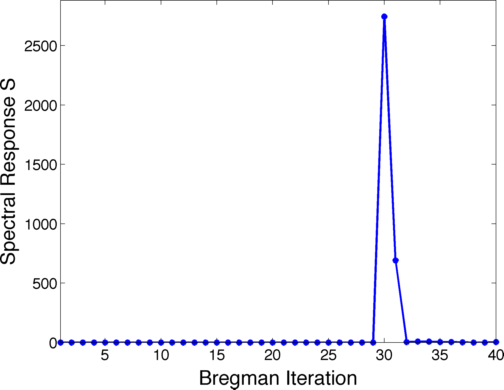}}\quad 
        \subfigure[Spectral response $S$ for (c).]{\includegraphics[width=0.2\textwidth]{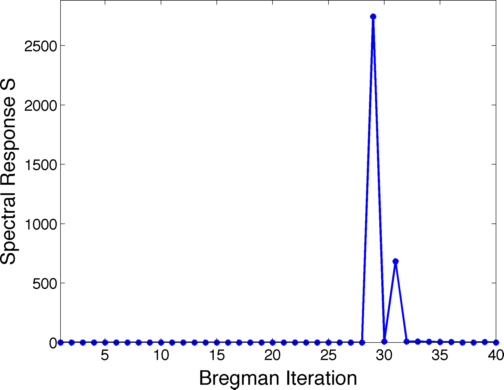}}
  \caption{\textit{Scale ambiguity w.r.t. size and intensity.} (a) - (c) Three different datasets that differ only in the intensity of the larger circle (0.68/0.69/0.7). The associated spectral responses $S$ are shown below. (a) together with (d) shows that out method is not able to reliably distinguish size and intensity scales but only a slight decrease of the intensity leads to two separate peaks in the spectral response function ((b) together with (e) and (c) together with (f)). }
\label{fig:criticalcase} 
\end{figure}
Figure \ref{fig:criticalcase}(a)-(c) shows three different data sets, each with a small disc with intensity 1 and a larger discs with intensity $0.68-0.7$ respectively. For all three examples we computed the spectral response with fixed parameters $\alpha = 200$ and 40 Bregman iterations shown in Figure \ref{fig:criticalcase}(d)-(f). We can see that in (d) the spectral method was not able to detect the difference between a large object with a lower intensity and a smaller object with a high intensity. Both objects appear at exactly the same point in time. Figure \ref{fig:criticalcase}(e) and (f) shows that already a small increase of the lower intensity circumvents the problem and both objects are represented as individual peaks in the spectral response. The height of each peaks reveals which object occurs. Yet we hope to solve this issue in future versions; solution ideas are presented in the summary and conclusion section \ref{sec:conl}.

\subsection{Robustness against Noise}
In this second subsection we investigate the robustness of our methods with respect to uncertainties like noise. 
\begin{figure}[htb]
  \centering 

  \subfigure[Noise level $\sigma = 0.25$.]{\includegraphics[width=0.22\textwidth]{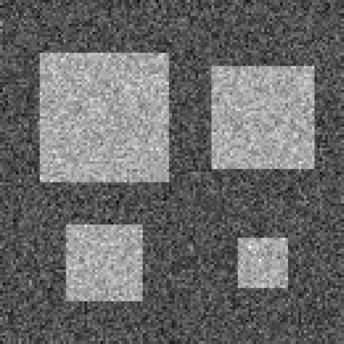}}\quad 
  \subfigure[Noise level $\sigma = 0.5$.]{\includegraphics[width=0.22\textwidth]{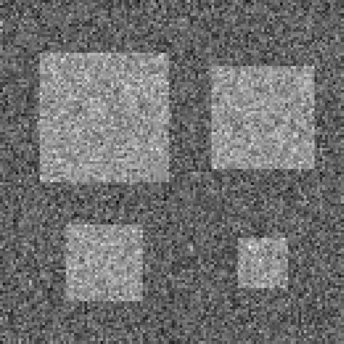}}\quad 
  \subfigure[Noise level $\sigma = 0.75$.]{\includegraphics[width=0.22\textwidth]{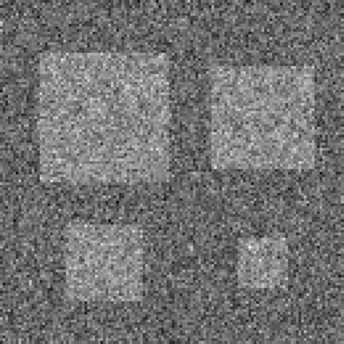}}\quad 
  \subfigure[Noise level $\sigma = 1$.]{\includegraphics[width=0.22\textwidth]{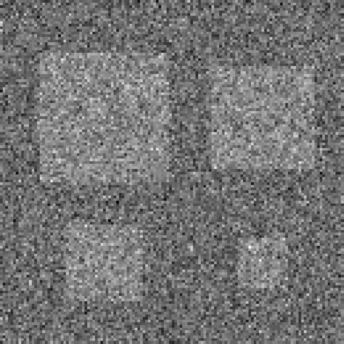}}\\
   \subfigure[Spectral response $S$.]{\includegraphics[width=0.22\textwidth]{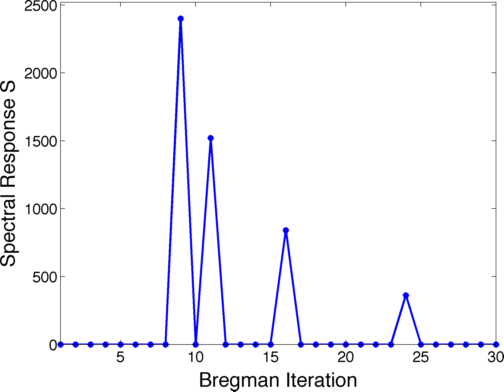}}\quad 
  \subfigure[Spectral response $S$.]{\includegraphics[width=0.22\textwidth]{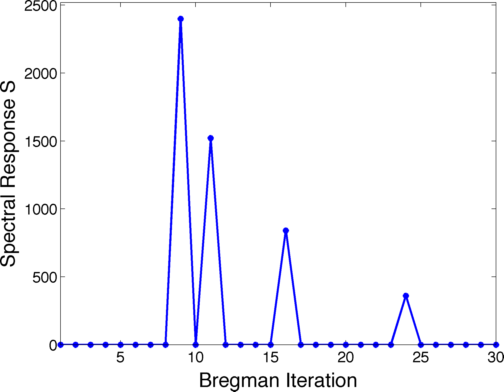}}\quad 
  \subfigure[Spectral response $S$.]{\includegraphics[width=0.22\textwidth]{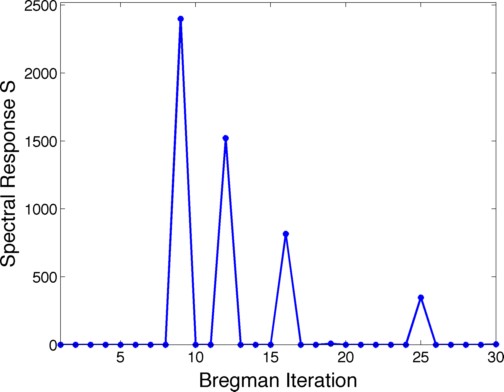}}\quad 
  \subfigure[Spectral response $S$.]{\includegraphics[width=0.22\textwidth]{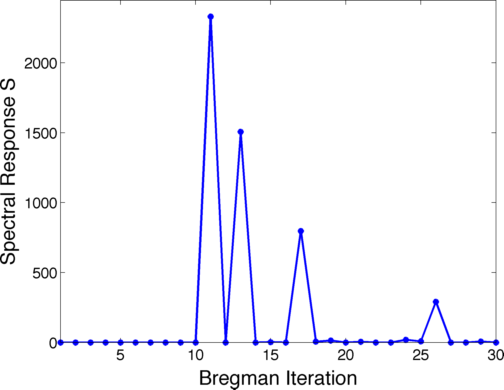}}\\
   \subfigure[Color-coded segmentation of multiple scales.]{\includegraphics[width=0.22\textwidth]{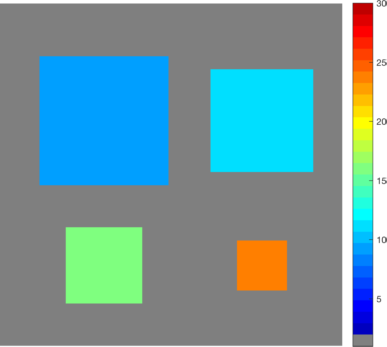}}\quad 
  \subfigure[Color-coded segmentation of multiple scales.]{\includegraphics[width=0.22\textwidth]{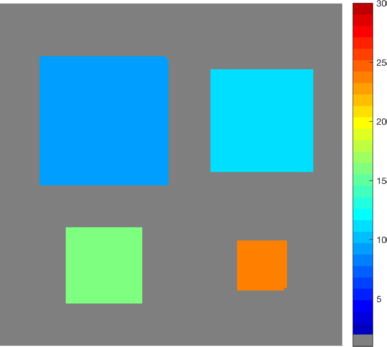}}\quad 
  \subfigure[Color-coded segmentation of multiple scales.]{\includegraphics[width=0.22\textwidth]{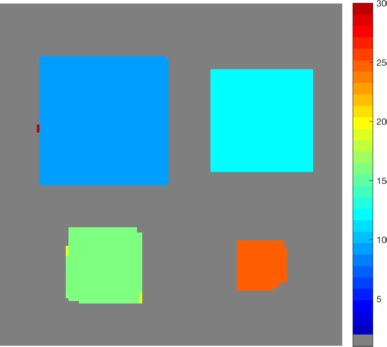}}\quad 
  \subfigure[Color-coded segmentation of multiple scales.]{\includegraphics[width=0.22\textwidth]{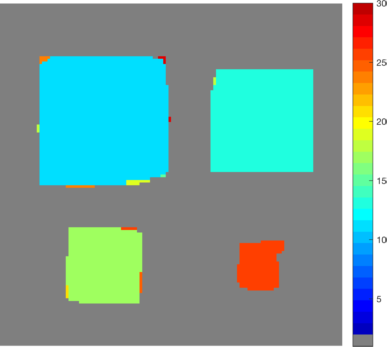}}
      
\caption{\textit{Robustness against noise}. (a) - (d) shows the same underlying binary signal but with different levels of noise added. Exactly the same model parameters are used to process all 4 datasets. The associated spectral response functions are shown in (e) - (h) and the color-coded segmentation results are shown in (i) - (l). We can see that our method is very robust against noise, we only see a slight shift to the right in the spectral response function for higher noise levels. For high noise levels some artifacts at the objects boundary occur but we reliably find all 4 objects independent of the noise level. }
\label{fig:robustness} 
\end{figure}
Figure \ref{fig:robustness}(a)-(d) shows the same underlying image but with different noise levels. The underlying signal is a binary signal showing four squares with different size. We added Gaussian noise with mean $0$ and standard deviation $\sigma = 0.25, 0.5, 0.75$ and $1$ respectively. Since squares are the unit balls of the $\ell^{\infty}$ norm we chose $\gamma = \| \cdot \|_{1}$ and therefore $\gamma^{\ast} = \| \cdot \|_{\infty}$. The results of our inverse Bregman-CV model are plotted in terms of the spectral response function $S$ and the color-coded segmentation of all scales $\int_{t}\Phi(x,t) \cdot t \mathrm{d} t$. Note that we chose the same parameter setting for all four examples ($\alpha = 100$ and 30 Bregman iterations). In Figure \ref{fig:robustness}(e)-(h) we can see that the spectral response function has for all noise levels a sparse representation indicating very clearly the four scales in the underlying signal. It is also remarkable that the peaks are roughly at the same points in time, shifting only slightly to the right for higher noise levels. This is even more obvious in the color-coded representation Figure \ref{fig:robustness}(i)-(l). The color bar next to each image shows that the color coding is consistent over all four examples. Thus, the same colors (see for example Figure \ref{fig:robustness}(i) and (j)) reflect that the objects occur at the exact same iteration unaffected of the changed noise level. Even in case of high noise ($\sigma$ = 0.75) only the second and the fourth largest object appear one iteration later although the smaller objects are now affected by an incomplete segmentation. In case of very high noise ($\sigma = 1$) we cannot perfectly segment the objects, especially the smaller ones "get lost" in the noise. Nevertheless the spectral response function is still close to sparse and for all four objects we reliably segment most of the object. These results show that our method is very robust to noise and circumvents to some extend the problem of choosing a good regularization parameter $\alpha$. The parameter setting that we used yielded a reliable extraction of all scales in the underlying signal independent of the noise present in our data. 
\subsection{Segmentation regarding Shapes of Eigenfunctions}\label{sec:resultsshapes}
This subsection is about the choice of the $\gamma$-function and how this choice is reflected in our spectral analysis framework. As mentioned in Section \ref{sec:genTV} by varying $\gamma$ in the generalized definition of the total variation \eqref{eq:TVgen} the shape of the eigenfunctions is changed. In the following two examples we will first show different examples of eigenfunctions that lead to a sparse representation of the spectral response function. Here we differentiate between eigenfunctions that are Wulff shapes of the chosen $\gamma$ function and those that are not. In addition, we investigate what happens if the choice of $\gamma$ does not match the shapes present in the image.

\begin{figure}[htb]
  \centering 
  \subfigure[Data $f$ with segmentation shown in red.]{\includegraphics[width=0.22\textwidth]{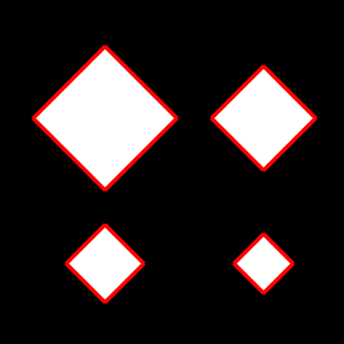}}\quad 
  \subfigure[Data $f$ with segmentation shown in red.]{\includegraphics[width=0.22\textwidth]{results_new/balls_size_nonoise/bregman_cv_24.png}}\quad 
  \subfigure[Data $f$ with segmentation shown in red.]{\includegraphics[width=0.22\textwidth]{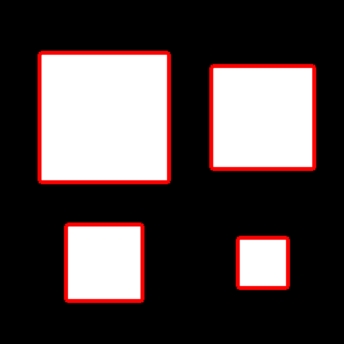}}\\
    \subfigure[Spectral response $S$.]{\includegraphics[width=0.22\textwidth]{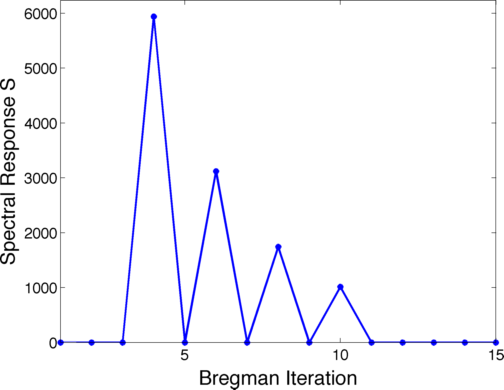}}\quad 
      \subfigure[Spectral response $S$.]{\includegraphics[width=0.22\textwidth]{results_new/balls_size_nonoise/S_gesamt.png}}\quad 
        \subfigure[Spectral response $S$.]{\includegraphics[width=0.22\textwidth]{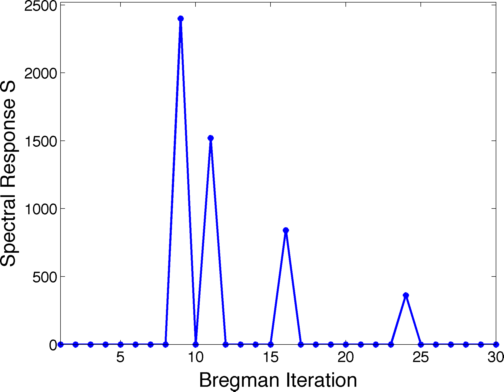}}\\
    \subfigure[Color-coded segmentation of multiple scales.]{\includegraphics[width=0.22\textwidth]{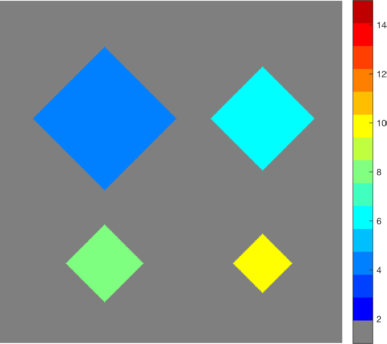}}\quad 
      \subfigure[Color-coded segmentation of multiple scales.]{\includegraphics[width=0.22\textwidth]{results_new/balls_size_nonoise_new/back_gray.png}}\quad 
        \subfigure[Color-coded segmentation of multiple scales.]{\includegraphics[width=0.22\textwidth]{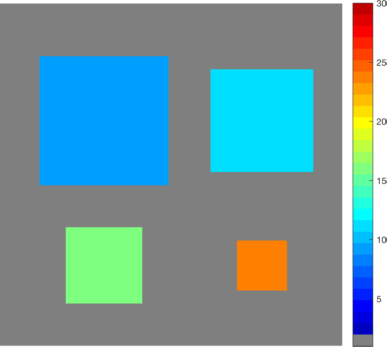}}
  \caption{\textit{Varying norm bodies and eigenshapes.} Automatic detection of different signals composed of eigenshapes for different norm bodies. Column one shows the results for $\gamma^{\ast} = \| \cdot \|_1$, column two for $\gamma^{\ast} = \| \cdot \|_2$ and column three for $\gamma^{\ast} = \| \cdot \|_{\infty}$. (a) - (c) Segmentation results visualized via red contour after 30 Bregman Iterations. Associated spectral response function $S$ is shown in (d) - (f) and the color-coded segmentation in (g) - (i). With the correct choice of $\gamma$ all eigenshapes can be detected in a single step. }
\label{fig:differentnorms} 
\end{figure}
In Figure \ref{fig:differentnorms}(a)-(c) each images shows a composite of eigenfunction for three different $\ell^p$-norm choices for $\gamma$ (namely $p = \infty$, $2$ and $1$). In each of the images, the input data is shown together with the final segmentation result in red. The associated spectral responses are shown below (Figure \ref{fig:differentnorms}(d)-(f)). We can see that as long as the choice of $\gamma$ and the given dataset match (in the sense that the components of $f$ have the same shape as the unit-ball of $\gamma^{\ast}$ (cf. Table \ref{tab:shapes})) we get a sparse response function. In Figure \ref{fig:differentnorms}(g)-(i) the color-coded segmentation of all scales is shown. We can see that for rectangular shapes ((g) and (i)) no discretization artifacts at the boundary occur. Moreover, we can again use the same parameter setting for all three data sets. The similar color coding in (h) and (i) indicates that the objects occur roughly at the same time point independent of their shape and $\gamma$. The data set in the first column is larger than the other ones thus the objects occur earlier and we need less iterations. 

\begin{figure}[t]
  \centering 
  \subfigure[TV eigenfunction for $\gamma = \| \cdot \|_2$.]{\includegraphics[height=0.14\textheight]{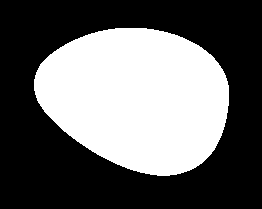}}\quad 
   \subfigure[Spectral response $S$.]{\includegraphics[height=0.14\textheight]{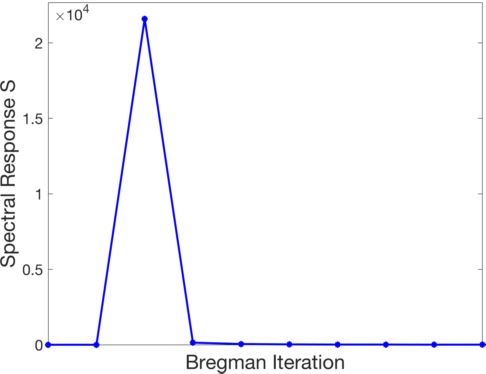}}\quad 
 \subfigure[Color-coded segmentation of multiple scales.]{\includegraphics[height=0.14\textheight]{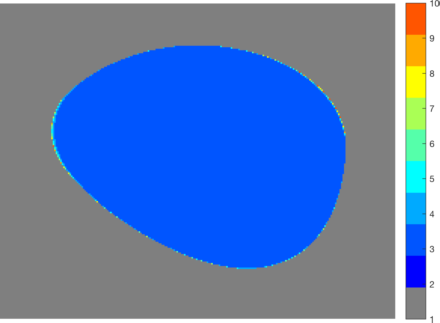}}\\
    \subfigure[TV eigenfunction for $\gamma = \| \cdot \|_1$.]{\includegraphics[width=0.235\textwidth]{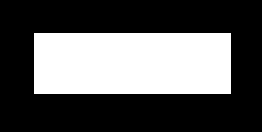}}\quad 
 \subfigure[Spectral response $S$.]{\includegraphics[height=0.14\textheight]{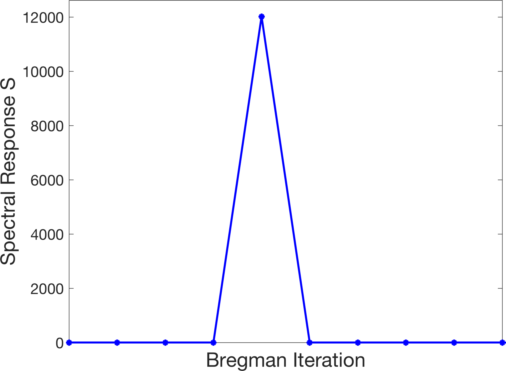}} \quad 
 \subfigure[Color-coded segmentation of multiple scales.]{\includegraphics[height=0.091\textheight]{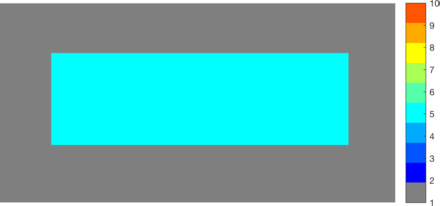}}
  \caption{\textit{Reconstruction of eigenfunctions that are not Wulff shapes.} (a) and (d) show two TV eigenfunctions that are not Wulff shapes for $\gamma = \| \cdot \|_2$ and $\gamma = \| \cdot \|_1$ respectively. The associated spectral responses are shown in (b) and (e) and the color-coded segmentation results in (c) and (f). It is obvious that eigenfunctions that are not Wulff shapes do have a sparse spectral response representation and appear in one step in the corresponding segmentation.}
\label{fig:nonWulff} 
\end{figure}
 In Figure \ref{fig:nonWulff} (a) and (d) we present two eigenfunctions which are not Wulff shapes. The round object in (a) fulfills condition \eqref{eq:l2eig} and is therefore a TV eigenfunction for $\gamma = \| \cdot \|_2$. The corresponding spectral response function in (b) and the color-coded segmentation in (c) show that not only eigenfunctions that are Wulff shapes appear in one step but also more general eigenfunctions. Another example (now for $\gamma = \| \cdot \|_1$) is shown in Figure \ref{fig:nonWulff} (d). Esedoglu and Osher showed in \cite{esedoglu2004} that a rectangle whose contours are parallel to the x- and y-axis respectively is an eigenfunction of the anisotropic TV functional. This is reflected in the spectral response function Figure \ref{fig:nonWulff} (e) and the color-coded segmentation in \ref{fig:nonWulff} (f) which show that the rectangle appears in one step although it is not a Wulff shape of the 1-norm. For natural images, where objects are often not Wulff shapes, this allows objects with various shapes to appear in one step. Some examples can be seen in the subsections \ref{sec:resultsCTC} - \ref{sec:resultsnetwork}.

\begin{figure}[t]
  \centering 
  \subfigure[Given Data.]{\includegraphics[width=0.23\textwidth]{results_new/varying_shapes_nonoise/bregman_cv_1.png}}\quad 
  \subfigure[Bregman iteration 10.]{\includegraphics[width=0.23\textwidth]{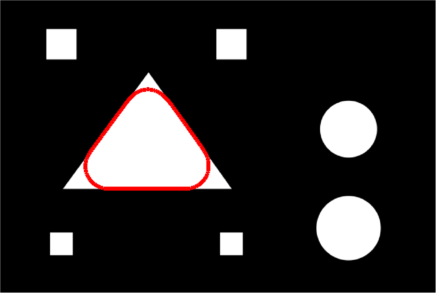}}\quad 
  \subfigure[Bregman iteration 16.]{\includegraphics[width=0.23\textwidth]{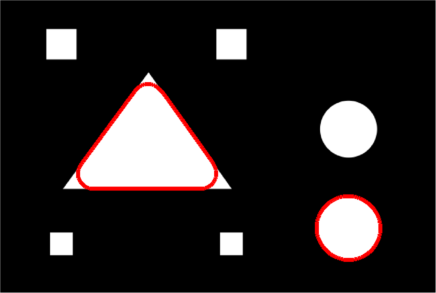}}\\
  \subfigure[Bregman iteration 18.]{\includegraphics[width=0.23\textwidth]{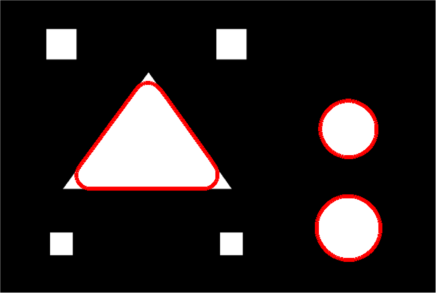}}\quad 
  \subfigure[Bregman iteration 31.]{\includegraphics[width=0.23\textwidth]{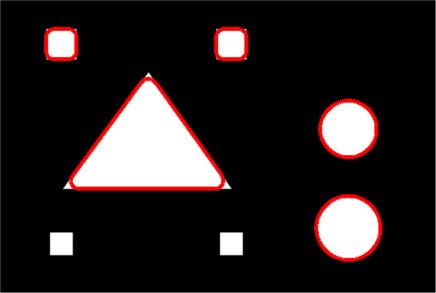}}\quad 
  \subfigure[Bregman iteration 41.]{\includegraphics[width=0.23\textwidth]{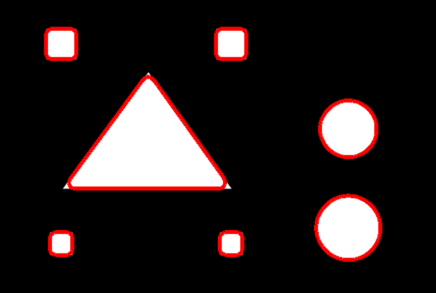}}\\
  \subfigure[Spectral response $S$.]{\includegraphics[width=0.23\textwidth]{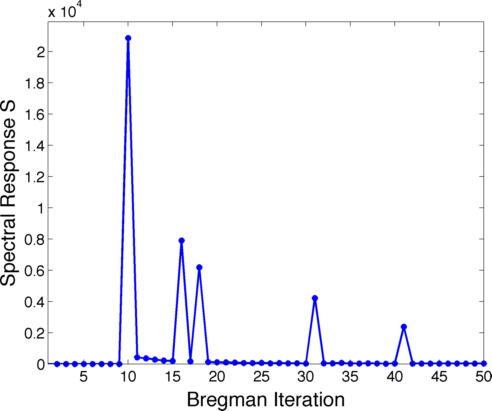}}\quad
  \subfigure[Color-coded segmentation of multiple scales.]{\includegraphics[width=0.31\textwidth]{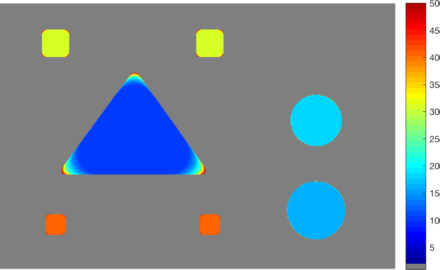}}  
  \caption{\textit{Reconstruction of mixed shapes.} Automatic detection of a signal consisting of different size scales and shapes. The binary signal is shown in (a). In (b) - (f) the segmentation results of the Bregman iterations where a significant peak in the spectral response $S$ in (g) appeared are shown. The color-coded segmentation is shown in (h). Here the behavior of non-eigenfunction over time can be seen. They do not appear in one step but reshape and evolve into the shape of the original signal.  }
\label{fig:scalesmix} 
\end{figure}
In Figure \ref{fig:scalesmix} we investigated how our multiscale approach behaves in case of shape mixtures and shapes that do not match the prior shape assumption made in $\gamma$. The given data (Figure \ref{fig:scalesmix}) consists of different size scales and different sizes with only two objects matching the prior assumption of round objects ($\gamma = \| \cdot \|_2$). Figure \ref{fig:scalesmix}(b)-(f) shows all time steps where a new scale occurs in the segmentation (shown in red). We can see that the round objects appear in one step while all other ones appear with rounded edges and reshape over time. This can be seen more clearly in the color-coded segmentation of all steps shown in Figure \ref{fig:scalesmix}(h). Especially for the large triangle in the middle we can nicely see how non-eigenfunctions are reshaped over time and how the segmentation is very round at the beginning and then propagates into the edges. Nevertheless we cannot get clear edges, neither for the triangle nor for the squares. The spectral transform function is nevertheless close to sparse but due to this reshaping behavior there is a small signal also between the clear peaks. It might be suitable to use the $\ell_0$-norm of $S$ as a measure in order to evaluate if a certain signal matches the prior shape assumption made. However for this idea the influence of discretization artifacts as shown in some examples before, would have to be clarified in further detail.
\subsection{Fluorescence microscopy images containing Circulating Tumor Cells}\label{sec:resultsCTC}
In this subsection we demonstrate how our multiscale segmentation approach can be applied to the analysis of Circulating Tumor Cells (CTCs). Before presenting some experimental data sets we will first give a short introduction on CTCs, associated research questions and challenges for diagnosis and treatment of cancer patients.

\begin{figure}[ht]
  \centering 
  \subfigure[Fluorescence image from the CellSearch system.]{\includegraphics[width=0.25\textwidth]{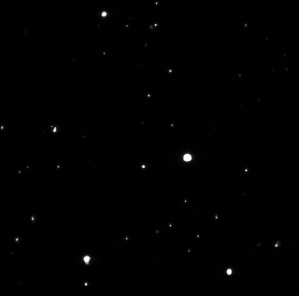}}\quad 
  \subfigure[Color-coded spectral response $S$.]{\includegraphics[width=0.25\textwidth]{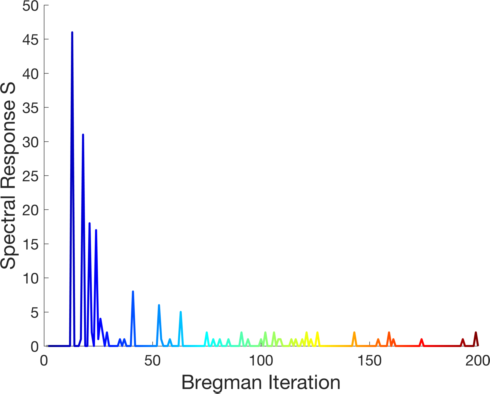}}\quad 
  \subfigure[Color-coded segmentation of multiple scales.]{\includegraphics[width=0.25\textwidth]{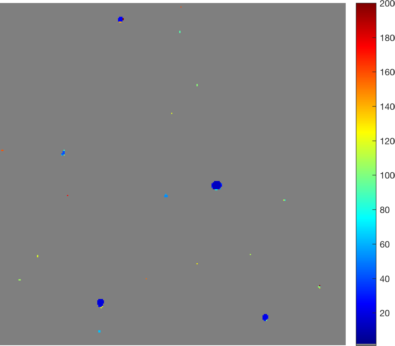}}\\
  
    \subfigure[Fluorescence image from a filter system.]{\includegraphics[width=0.25\textwidth]{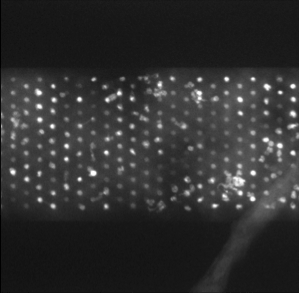}}\quad 
  \subfigure[Color-coded spectral response $S$.]{\includegraphics[width=0.25\textwidth]{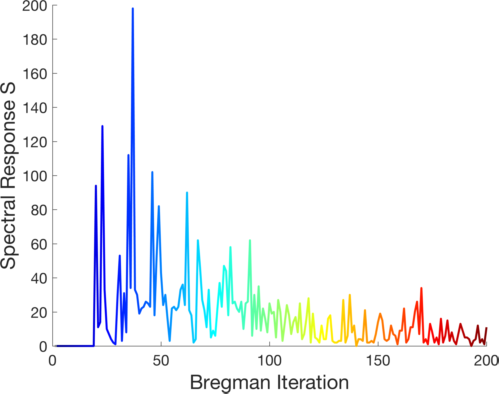}}\quad 
  \subfigure[Color-coded segmentation of multiple scales.]{\includegraphics[width=0.25\textwidth]{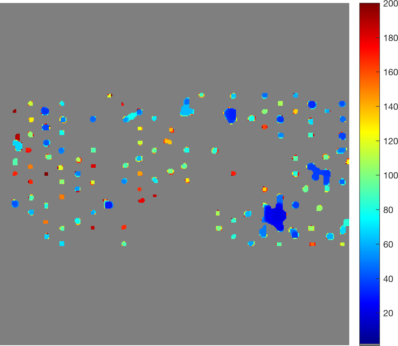}}\\
  
    \subfigure[Fluorescence image from the CellSearch system with more clustered cells.]{\includegraphics[width=0.27\textwidth]{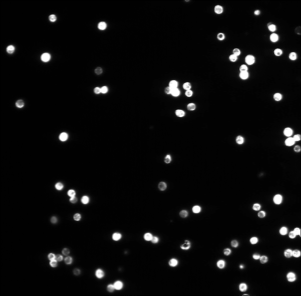}}\quad 
  \subfigure[Color-coded spectral response $S$.]{\includegraphics[width=0.25\textwidth]{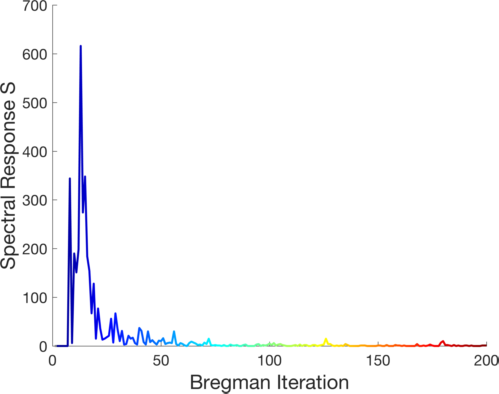}}\quad 
  \subfigure[Color-coded segmentation of multiple scales.]{\includegraphics[width=0.25\textwidth]{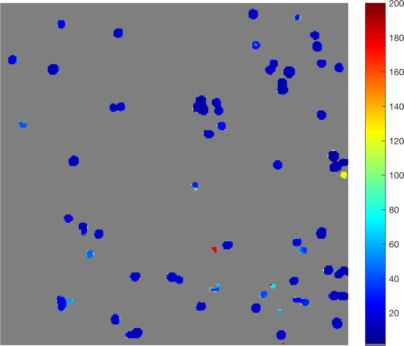}}\\
  
    \subfigure[Fluorescence image from the ARIOL system.]{\includegraphics[width=0.25\textwidth]{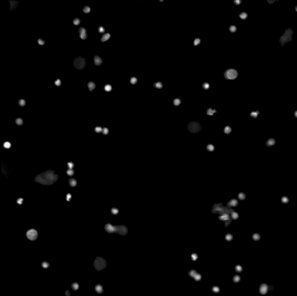}}\quad 
  \subfigure[Color-coded spectral response $S$.]{\includegraphics[width=0.25\textwidth]{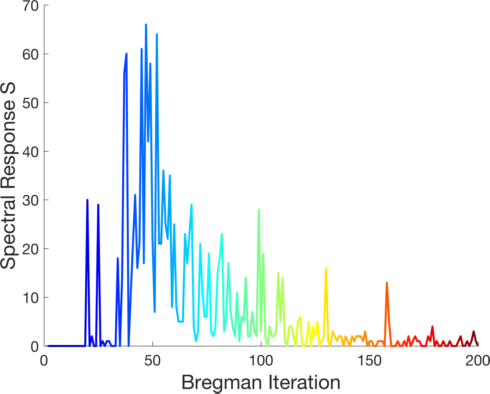}}\quad 
  \subfigure[Color-coded segmentation of multiple scales.]{\includegraphics[width=0.25\textwidth]{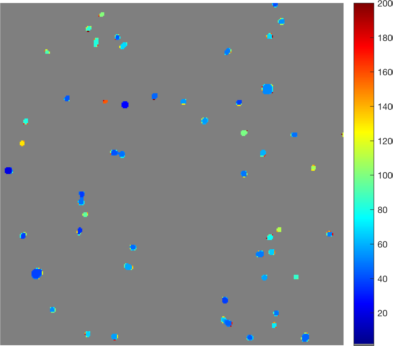}}\\

   \caption{\textit{Experimental cell data.} Automatic detection of (a) different sparse cell sizes (scales) with nearly constant background, (d) cells in front of an inhomogeneous background, (g) clustered cells with nearly constant background and (j) different intensity and size scales. All results are obtained with the exact same parameters. In column one single channel images obtained from different fluorescence microscopes are shown. The corresponding spectral responses are shown in the second column and the color-coded segmentation in column three. It can be seen that more complex images (inhomogeneous background or mixtures of scales) lead to more complexity in the spectral response function (response becomes less sparse).}
\label{fig:cells} 
\end{figure}

Circulating Tumor Cells are cells that dissociate from a primary tumor and invade the blood stream. In recent literature it was shown that the number of CTCs present in the blood is associated with the survival chance of a patient and can be used to guide therapy of cancer patients. The elimination of CTCs indicates effective therapy whereas increase or failure to eliminate indicates a futile therapy. A challenge in the identification of CTCs is that they are very rare and therefore difficult to detect among other cells in the pool of blood cells. The gold standard for CTC enumeration is the CellSearch system and in prospective multicenter studies a threshold of 5 CTCs / 7.5 ml of blood was used to separate patients with metastatic disease into those with favorable and unfavorable prognosis \cite{cristofanilli2004,debono2008}. Currently all known CTC analysis tools are based on subjective morphological criteria and objects can be classified differently by different operators which results in different interpretation of the status of the patients. The low cell counts and subsequent statistical analyses further increases the chances of mistakes. The problem of subjectivity in the CTC analysis is currently addressed by the development of an open-source toolbox to automatically detect and classify CTCs in datasets from various machines and institutes. This project is part of a European consortium called CANCER-ID that aims at validating blood-based biomarkers and is funded by the European Union. One key component of this toolbox is the development of a (nearly) parameter-free segmentation algorithm which performs reliably and efficiently on multiple tumor cell datasets (see Figure \ref{fig:cells}).
\begin{figure}[t]
\centering
 \includegraphics[width=0.4\textwidth]{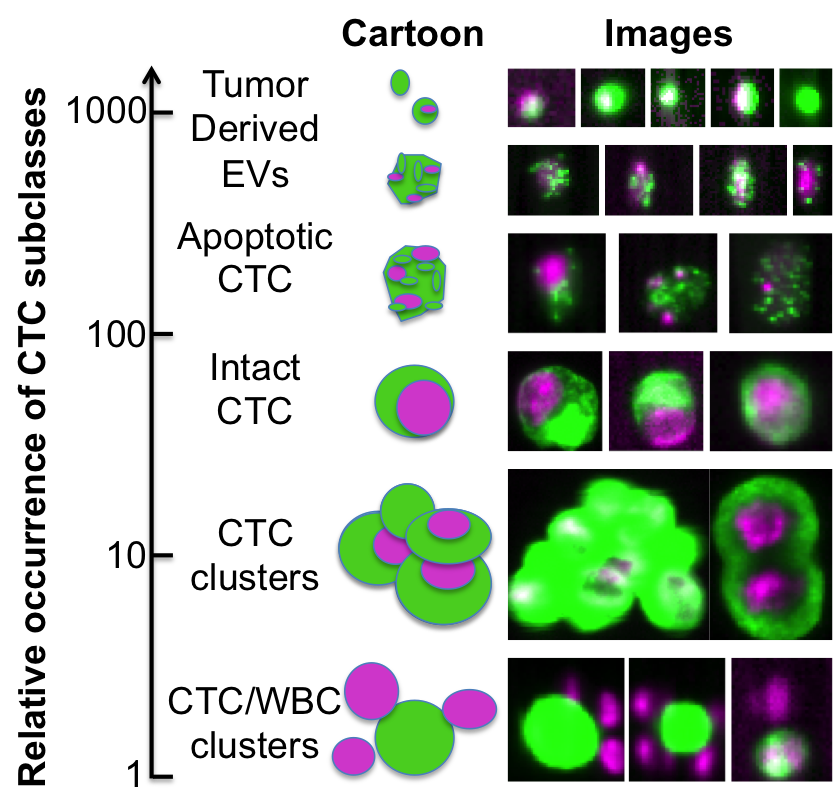}
   \caption{\textit{Schematic overview of different CTC classes.} For each subclass a cartoon image of a typical cell is shown. Next to it on the right different examples for every class are shown. All images consist of two overlaid color channels where green is the tumor cell marker and magenta is the nucleus marker. The relative frequencies of each subclass are shown on the left. By looking only at the green signal it is obvious that the size scale varies between different subclasses.}
\label{fig:cellclasses} 
\end{figure}
Moreover, in recent studies \cite{coumans2010} it was hypothesized that not only the definition used to classify objects as CTC predicted clinical outcome but also those objects that did not match the criteria of an intact cell. A schematic overview of different cell classes and their relative frequencies is shown in Figure \ref{fig:cellclasses}. Here, green is a fluorescence marker indicating the presence of cytokeratin expressed in cancer cells but not on white blood cells and magenta is a fluorescence marker indicating the presence of DNA thus showing the nucleus of the cell. When comparing the green marker among the different CTC classes we observe that the signal clusters strongly vary in size (area) and therefore motivate the usage of our multiscale segmentation approach.

In Figure \ref{fig:cells} (a), (d), (g) and (j) experimental data sets are shown. From the original image sets, consisting of three different fluorescent color channels, we extracted the tumor cell marker (green) and used those images as input for our multiscale segmentation approach. Although the difficulties vary between all images (inhomogeneous background, noise, cell clusters, mixture of size and intensity scales), we can process all images with our multiscale segmentation approach with the \textit{exact same parameters}. This is essential  for the development of a user-friendly (parameter-free) toolbox for CTC analysis. Note that the dim spots in image (d) are not cells but only pores of a filter used to collect cells (bright spots) and therefore it is not desired to segment them. The resulting spectral response functions for all four images are shown in Figure \ref{fig:cells} (b), (e), (h) and (k) with a color coding corresponding to the coding used in the segmentation results in (c), (f), (i) and (l). The color coding of the response function shows that all objects which appear later in our segmentation and therefore belong to finer scales have a yellow to reddish color in the color-coded segmentation. The very large and intact cells are blue (with some small artifacts at the boundary) and smaller cells (or large fragments) are shown in light blue to green. We can nicely observe that the object colors cover the whole color scale range. For images that are more complex (e.g. (d) and (j)) also the spectral response function is more complex but the color-coded segmentation shows that nearly every object appears in one step and thereby has a clearly defined scale that we use as a feature in our classification approach. Here, we profit from the fact that not only Wulff shapes (perfectly circular objects) but also other eigenshapes appear in one step in our segmentation.
Hence, this segmentation approach not only provides all contours automatically without any parameter adaptations but simultaneously also a simple classification of cells based on their size (scale resp. color/appearance time). This analysis can be applied to all color channels separately and be used together with more features in a subsequent automatic classification approach. The constants $c_1$ and $c_2$ can again simply be estimated a-priori from the data by a simple thresholding and averaging approach and are fixed throughout the iterative process.

\subsection{Electron microscopy cell images with nonuniform background}\label{sec:resultsEM}
\begin{figure}[t]
  \centering 
  \subfigure[Electron microscopy data of acinar cells in the parotid gland provided by \cite{riedel}.]{\includegraphics[height=0.15\textheight]{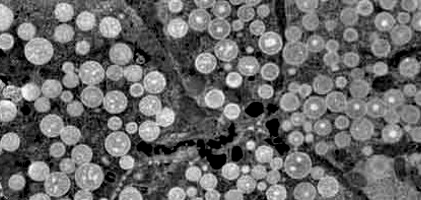}}\quad 
   \subfigure[Color-coded spectral response $S$.]{\includegraphics[height=0.15\textheight]{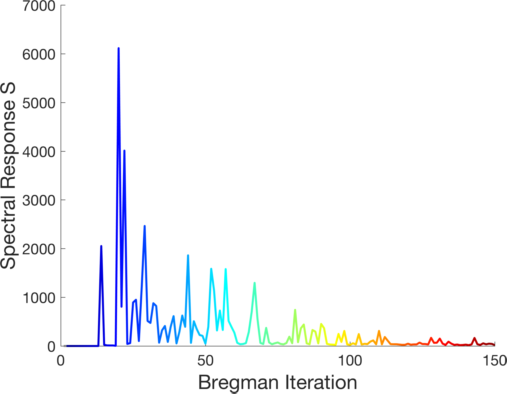}}\quad 
 \subfigure[Color-coded segmentation of multiple scales.]{\includegraphics[height=0.15\textheight]{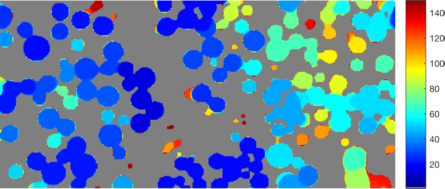}}
  \caption{\textit{Segmentation of cells in front of an inhomogeneous background.} (a) shows the input image obtained by electron microscopy, (b) the corresponding spectral response function and (c) the color-coded segmentation. This result shows that our method can detect inhomogeneities in the object we want to segment while inhomogeneities in the background (up to a certain degree) do not affect the result of our segmentation method. }
\label{fig:emcells} 
\end{figure}

In this subsection we show that our multiscale segmentation approach can also be used for real-world images where the background is less homogenous as in the cell images before. In Figure \ref{fig:emcells} (a) we show an electron microcopy dataset of acinar cells provided by \cite{riedel}. For a more obvious visual interpretation we reverted the contrast from the original EM dataset (bright background with dark cells) to a dark, but not uniform, background with brighter cells. The intensity of the cells is also inhomogeneous and varies from very bright cells (left) to cells with an intensity closer to the background (right). Figure \ref{fig:emcells} (b) shows the color-coded spectral response function and (c) the color-coded segmentation of the cells. We see that the inhomogeneous intensity of the cells is reflected in the scales found by our method. The inhomogeneity of the background is not influencing the segmentation result and the constants $c_1$ and $c_2$ can still be estimated beforehand from the data without further adaption. One problem that occurs is the ambiguity of size and intensity scales which can be observed when comparing small but bright cells on the right (visualized in light green) with larger cell clusters on the right (also visualized in light green). Currently our method is not able to differentiate their scales. A useful extension of our method would be a combination of our multiscale approach with a watershed type algorithm that is able to split cell clusters into individual cells of different scales.

\subsection{Segmentation and clustering of biological network structures}\label{sec:resultsnetwork}
\begin{figure}[t]
  \centering 
  \subfigure[Given Data.]{\includegraphics[height=0.2\textheight]{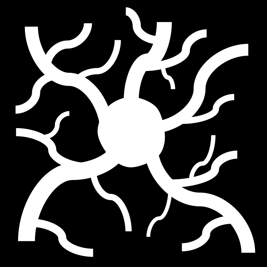}}\quad 
   \subfigure[Color-coded spectral response $S$.]{\includegraphics[height=0.2\textheight]{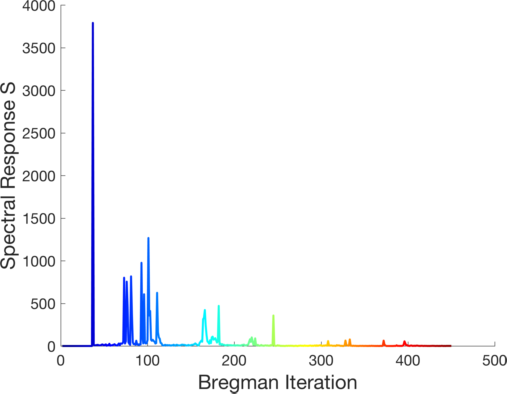}}\quad 
 \subfigure[Color-coded segmentation of multiple scales.]{\includegraphics[height=0.2\textheight]{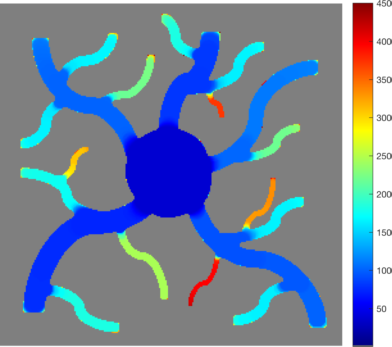}}
  \caption{\textit{Reconstruction of different scales in network-like structures.} The input data is shown in (a) with the corresponding spectral response function (b) and the color-coded segmentation (c). In this example the usability of our method in case of more complex, network-like structures is demonstrated. Interestingly, the spectral response function shows distinct sparse clusters where every cluster corresponds to "arms" or "branches" of a specific diameter. Although they are not shaped as an eigenfunction they appear (nearly) in one step and are represented by one specific scale.}
\label{fig:cellarms} 
\end{figure}

\begin{figure}[t]
  \centering 
  \subfigure[Retinal image of the STARE dataset \cite{hoover2000}.]{\includegraphics[height=0.2\textheight]{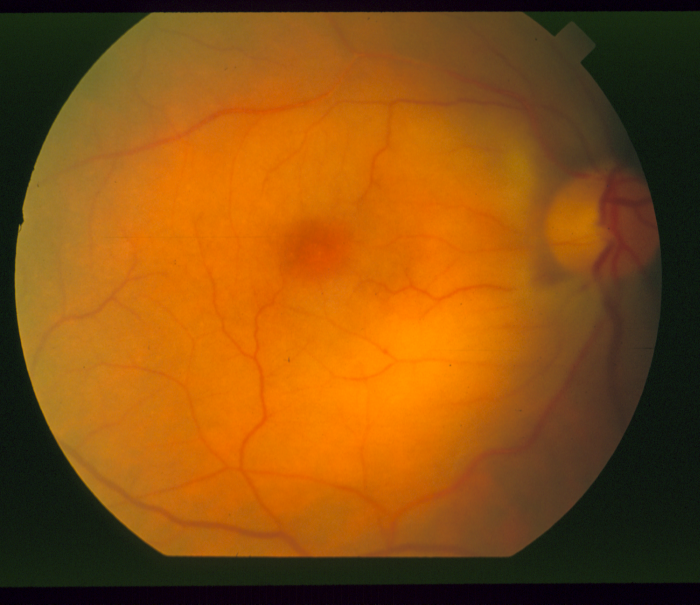}}\quad 
     \subfigure[Manually segmented blood vessel network \cite{hoover2000}.]{\includegraphics[height=0.2\textheight]{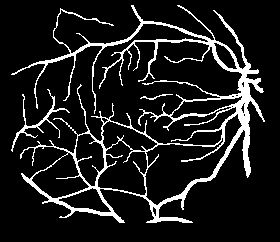}}\quad 
   \subfigure[Color-coded spectral response $S$.]{\includegraphics[height=0.2\textheight]{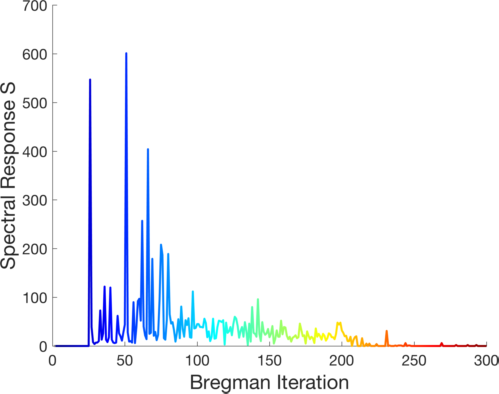}}\quad 
 \subfigure[Color-coded clustering of network in (b).]{\includegraphics[height=0.2\textheight]{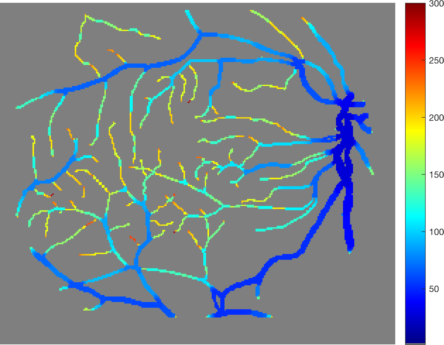}}
  \caption{\textit{Clustering of networks.}  (a) A retinal image of a human eye taken from the STARE dataset \cite{hoover2000} with a manual segmentation of the blood vessel network in (b). We applied our multiscale segmentation method to the binary image in (b) and thereby obtain the spectral response function (c) and a clustering of the network based on the diameter of the network branches (d). We see that large blood vessels are represented by blue colors, medium vessels by light blue to green colors and the very fine branches by yellow to red colors. Although the representation of (c) is not as sparse as in the previous example we can get a meaningful clustering of the graph based on the resulting colors in (d).}
\label{fig:retina} 
\end{figure}

In many biological and medical applications network-like structures occur as part of vascular systems. Examples are blood-vessels in retinal images, anisotropic structures in brain or myocardium imaging or at the microscopic level even cells consisting of a round core with several meandering arms of different diameters. Although these networks have a much more complex shape than the examples before, our multiscale method is very useful for a segmentation or clustering of these networks. An example of a round object with "arms" of different diameters is shown in Figure \ref{fig:cellarms} (a). In (b) the corresponding spectral response function is shown and the color-coded segmentation in (c). We can see that the spectral response function has sparse clusters where each cluster corresponds to branches of a certain diameter. Although these "arms" are not eigenfunctions, they do appear in one step and therefore have a distinct scale if they have constant diameter. A similar behavior can be seen in Figure \ref{fig:retina}. As input for our multiscale method we used a manual blood vessel segmentation (b) of a retinal image (a) from the STARE dataset \cite{hoover2000}. The resulting spectral response is shown in Figure \ref{fig:retina} (c) and the color-coded segmentation in (d). We see that the different times of appearance of the vessels (indicated by the different colors) lead to a clustering of the underlying network based on the diameter of each branch. With this approach we can get segmentations of the network that contain only branches of certain diameter scales or directly get an estimate of how many different vessel sizes are present and how often they occur. In order to apply our segmentation approach directly on raw retinal images, instead of using it as a postprocessing step, a decomposition method for retinal images similar to Zosso \cite{zosso2015} could be combined with our multiscale approach.

\section{Summary and Conclusion}\label{sec:conl}

In this paper, we have studied a novel multiscale segmentation method which is a combination of an inverse scale space variational approach and a nonlinear spectral analysis. Using our approach we can automatically detect objects of different scales where scale can refer to intensity or size. These scales can be evaluated using a spectral response function and we can easily decompose the segmentation with respect to those scales using the spectral transform function. Thus we have extended the useful spectral analysis that was introduced for $TV$ denoising applications before to a segmentation approach. This novel approach can to some extent circumvent the problem of choosing a suitable regularization parameter. The iterative procedure using Bregman distances can be interpreted as an automatic regularization parameter choice method and is very robust with respect to noise or errors in the estimated $c_1$ and $c_2$ values in the CV model. We have shown that our method works very reliable even for applications with more than two intensity scales although this actually contradicts the models binary assumption. The framework can easily be adapted to different prior shape assumptions (Wulff shapes). With regard to real experimental datasets such an assumption is not limiting the usefulness of the method, due to its ability to even detect compositions of eigenfunctions not belonging to the restricted class of Wulff shapes. For the numerical implementation we used a primal-dual optimization scheme that, especially with a C/MEX version, lead to a very efficient framework, even for large datasets. We illustrated the strengths and limitations of our method on synthetic datasets with a certain focus on eigenfunctions, shapes as well as robustness against noise and background artifacts. Three different biomedical imaging applications with different shapes and challenges emphasize the potential and wide applicability of our approach. The code is also part of an open-source toolbox called ACCEPT. With this software users are able to automatically detect and classify CTCs from fluorescence images, but can also extract important information related to the level of expression of certain therapy targets on the CTC (see the application and results described in section \ref{sec:resultsCTC}).

One limitation that we currently have in our model is an ambiguity in scale. It is known for variational methods with an $l_2$ type data fidelity term combined with total variation regularization that these methods are not able to distinguish between intensity and size scales. In practice it seems not very likely that two objects with different size and different intensity fall onto exactly the same peak in the spectral response. Nevertheless to reliably interpret the spectral response function it is advantageous to have only one dominating scale in the data that is analyzed. Another limitation is the flexibility with regard to different shapes. Currently only one shape assumption at the moment is possible.

\textit{Outlook.} There are four main questions that we want to address in our future research. To allow a clear interpretation of the spectral response function even in the case of size and intensity scales we want to adapt the dataterm in the way that it can differentiate both scales. For denoising this can be done using an $l_1$ type data fidelity and we want to generalize this to our segmentation model (which is not obvious in terms of gradient flows). For the experimental cell datasets we currently analyze every fluorescent channel separately. In a future work we would like to investigate if a combined segmentation in all channels can lead to an improved segmentation and a richer response function. Moreover, we are interested in how far we can use the spectral response function to investigate if a shape assumption matches the data. One idea is to use the $l_0$ norm of the spectral response function as a measure for shape optimization. A further development of this idea could lead to an approach where we could learn an optimal $\gamma$ function and thereby an optimal shape of eigenfunctions for a specific dataset. Another very interesting consideration is to transfer our method from the currently local setup to a nonlocal setup. A first foundation for this was recently laid with a work on non-local TV spectral theory by Aujol et al. in \cite{Aujol2015}. More mathematical spectral theory for nonlocal TV, its usefulness for solving and improving complex segmentation tasks and a direct comparison with spectral clustering and graph cuts will be of high interest for the community and our future research.
%
%
\section*{Acknowledgements}%
LZ, LT, CB acknowledge support by the EUFP7  program \#305341 CTCTrap and the IMI EU program \#115749 CANCER-ID. CB acknowledges support by the NWO via Veni grant 613.009.032 within the NDNS+ cluster.
%
%

\bibliographystyle{alpha} 
\bibliography{refs}

\end{document}

%% file: packages.tex
\usepackage{fullpage} 

\usepackage[affil-it]{authblk}%

\usepackage[pagebackref,pdftex,bookmarks,colorlinks,breaklinks]{hyperref}
\hypersetup{%
		linkcolor=magenta,      
		citecolor=blue,				
		filecolor=dullmagenta,
		urlcolor=magenta,
}

\usepackage[T1]{fontenc}    
\usepackage[utf8]{inputenc} 
\usepackage{lmodern}
\usepackage{ae}             

\usepackage{amsmath,amstext,amsfonts,amssymb,amsthm,latexsym}

\usepackage{todonotes}

\usepackage{graphicx}         
\usepackage{epsfig}           
\usepackage{subfigure}

\usepackage{algorithm}
\usepackage{algpseudocode}

\usepackage{bbm}



%% file: commands.tex
\newcommand{\eq}[1]{\begin{equation*}#1\end{equation*}}
\newcommand{\eqn}[1]{\begin{equation}#1\end{equation}}

\newcommand{\R}{\mathbb{R}}

\newcommand{\W}{\mathcal{W}}

\DeclareMathOperator*{\argmin}{\textnormal{argmin}}
\DeclareMathOperator*{\esssup}{\textnormal{ess~sup}}

\renewcommand{\div}{\nabla \cdot}%
\newcommand{\grad}{\nabla}%